\documentclass[twocolumn]{autart}
\usepackage[sort&compress]{natbib}

\usepackage{amsmath,amssymb,dsfont}
\usepackage[usenames]{color}
\usepackage{xcolor}
\usepackage{caption,subcaption}
\usepackage{graphicx}
\usepackage{bbm}
\usepackage{enumerate}

\usepackage{hyperref}
\hypersetup{
    colorlinks = true,
    allbordercolors = {white},
    allcolors = {purple},
}

\newcommand{\one}{{\mathbf{1}}}
\newcommand{\est}{{\hfill $\star$}}

\newcommand{\erem}{{\hfill $\lozenge$}}

\newcommand{\tran}{^{\top}}

\newcommand{\beq}{\begin{equation}}
\newcommand{\eeq}{\end{equation}}
\newcommand{\beqn}{\begin{equation*}}
\newcommand{\eeqn}{\end{equation*}}
\newcommand{\bea}{\begin{eqnarray}}
\newcommand{\eea}{\end{eqnarray}}
\newcommand{\beas}{\begin{eqnarray*}}
\newcommand{\eeas}{\end{eqnarray*}}
\newcommand{\ba}{\begin{array}}
\newcommand{\ea}{\end{array}}
\newcommand{\bit}{\begin{itemize}}
\newcommand{\eit}{\end{itemize}}
\newcommand{\ben}{\begin{enumerate}}
\newcommand{\een}{\end{enumerate}}
\newcommand{\dss}{\displaystyle}

\newcommand{\ped}[1]{ _{ {\mathrm{#1} } }}

\newcommand{\ap}[1]{ ^{ {\mathrm{#1} } }}

\newcommand{\Real}[1]{ { {\mathbb R}^{#1} } }
\newcommand{\Realp}[1]{ { {\mathbb R}_+^{#1} } }

\newtheorem{corollary}{Corollary}
\newtheorem{theorem}{Theorem}
\newtheorem{proposition}{Proposition}
\newtheorem{remark}{Remark}
\newtheorem{lemma}{Lemma}

\newcommand{\calV}{{\mathcal V}}
\newcommand{\calG}{{\mathcal G}}
\newcommand{\calE}{{\mathcal E}}
\newcommand{\calP}{{\mathcal P}}

\newcommand{\calT}{{\mathcal T}}

\newcommand{\ve}{\varepsilon}

\begin{document}
\begin{frontmatter}
\title{Clearing Payments in Dynamic Financial Networks\thanksref{footnoteinfo}} %

\thanks[footnoteinfo]{This paper was not presented at any IFAC meeting. }

\author[DET]{Giuseppe C. Calafiore}\ead{giuseppe.calafiore@polito.it},
\author[DET]{Giulia Fracastoro}\ead{giulia.fracastoro@polito.it},
\author[DET]{Anton V. Proskurnikov}\ead{anton.p.1982@ieee.org},

\address[DET]{Department of Electronics and Telecommunications, Polytechnic of Turin, Turin, Italy}  

\begin{abstract}
This paper proposes a novel dynamical model for determining clearing payments in  financial networks. We extend the classical Eisenberg-Noe model of financial contagion to multiple time periods, allowing financial operations to continue after possible initial {\em pseudo defaults}, thus permitting nodes to recover and eventually  fulfil their liabilities.
Optimal clearing payments in our model are computed by solving a suitable linear program, both in the full matrix payments case and in the pro-rata constrained case. We prove that the proposed model obeys the \emph{priority of debt claims} requirement, that is, each node at every step either pays its liabilities in full, or it pays out all its balance. In the pro-rata case, the optimal dynamic clearing payments are unique, and can be determined via a time-decoupled sequential optimization approach.

\end{abstract}
\end{frontmatter}


\section{Introduction}
The current global financial system
is  a highly interconnected network
of institutions that are linked together via a structure of mutual
debts or liabilities.
Such interconnected structure makes the system potentially prone
 to ``cascading defaults,'' whereby  a shock at a node (e.g., an expected incoming payment that gets cancelled or delayed for some reason) may provoke a default at that node, which then cannot pay its liabilities to neighbouring nodes, which in turn default, and so on in an avalanche fashion.
The global financial crisis of 2008 is an example of such behavior, where the bankruptcy of Lehman-Brothers is identified as the watershed event that started the crisis.
Since the consequences of these cascading events can be catastrophic, modeling and analyzing such behavior is of crucial importance. The seminal work~\cite{EisNoe:01} introduced a simple model for studying  financial contagion.   In particular, they focused on defining a clearing procedure between financial entities. Clearing consists in
a procedure for settling claims in the case of defaults, on the basis
of a set of rules and prevailing regulations.
 In~\cite{EisNoe:01}, the authors showed that there  exist a clearing vector which defines the mutual interbank payments, under certain assumptions.
 Among such assumptions, an important one is  that the debts of all nodes of the system are paid simultaneously.

The basic model presented in \cite{EisNoe:01} has become a cornerstone in the analysis of financial contagion and it has been extended in various directions. In particular,  non-trivial features were added in order to make the model more realistic. The models presented in  \cite{cifuentes2005liquidity,shin2008risk}, for instance, consider also the liquidity risk. Instead, in \cite{elsinger2009financial,suzuki2002valuing} cross-holdings and seniority of liabilities are introduced. Other works take into account costs of default \citep{rogers2013failure}, illiquid assets \citep{amini2016fully}, mandatory disclosures \citep{alvarez2015mandatory}, cross-ownership of equities and liabilities \citep{fischer2014no}, and decentralized clearing processes \citep{csoka2018decentralized}.

The vast majority of the works based on the Eisenberg-Noe model, however, considers the problem only in a static, or single-period, setting. This assumption is quite unrealistic, since it supposes that all liabilities are claimed and due at the same time. In addition, static models are only able to capture the immediate consequences of a financial shock.
For these reasons, several works recently  proposed time-dynamic extensions of the Eisenberg-Noe model. In \cite{sonin2017banks} a continuous-time model of clearing in financial networks is presented. This work has later been extended by considering liquid assets \citep{chen2021financial}, heterogeneous network structures over time and early defaults \citep{banerjee2018dynamic}. Other works \citep{feinstein2021dynamic} propose to combine the  interbank Eisenberg-Noe model and the dynamic mean field approach. Instead, \cite{feinstein2020capital} uses a continuous-time model for price-mediated contagion.

A different line of research extended the Eisenberg-Noe model considering a discrete-time setting. In \cite{capponi2015systemic,ferrara2019systemic} a multi-period clearing framework is introduced. Using a similar approach, \cite{kusnetsov2019interbank} considers the case where interbank liabilities can have multiple maturities, considering both long-term and short-term liabilities.

In the present work, we focus on a discrete-time setting and introduce a multi-period model whereby financial operations are allowed for a given number of time periods after the initial theoretical default (named here {\em pseudo default}).
This allows to reduce the effects of a financial shock, since some nodes may possibly recover and eventually fulfil their debts. We first consider the  general case
where
payment matrices are unconstrained.
%
%
This scenario has been introduced in the static case in \cite{CDC2021, Journal2021}, where its advantages over the proportional rule in terms of the overall system loss have been highlighted.
%
Here, we prove in a dynamic setting  that the optimal sequence of payment matrices satisfies the absolute priority of debt claims rule, hence the proposed method produces proper clearing matrices at each stage.

We then consider the situation in which a proportionality rule is enforced,
whereby nodes must  pay the claimant institutions proportionally to their nominal claims (pro-rata rule).
We prove that under the pro-rata rule the optimal payments are again proper clearing payments,
they are unique and, further, the multi-stage optimization problem can be decoupled in time into an equivalent series of LP problems.

The remainder of the paper is organized as follows. Section 2 introduces some preliminary notions and the notation that will be used in the next sections. In Section 3 we introduce the Eisenberg-Noe financial network model. Then, in Section 4 we illustrate the proposed dynamic model, considering both the unrestricted case and
the case with the pro-rata rule imposed. A schematic example is proposed in Section 5 in order to illustrate  the proposed model. Conclusions are drawn in Section 6. For ease of reading, we collected the proofs of all technical results in an appendix.

\section{Preliminaries and notation}

Given a finite set $\calV$, the symbol $|\calV|$ stands for its cardinality. The set of families $(a_{\xi})_{\xi\in\Xi}$, $a_{\xi}\in\mathbb{R}$, is denoted by $\Real{\Xi}$.
For two such families $(a_{\xi}),(b_{\xi})$, we write  $a\leq b$ ($b$ dominates $a$, or $a$ is dominated by $b$) if $a_{\xi}\leq b_{\xi}$,
$\forall\xi\in\Xi$.
 We write $a\lneq b$ if $a\leq b$ and $a\ne b$. The operations $\min,\max$ are also defined elementwise, e.g., $\min(a,b)\doteq(\min(a_{\xi},b_{\xi}))_{\xi\in\Xi}$.
 These notation symbols apply to both vectors (usually, $\Xi=\{1,\ldots,n\}$) and matrices (usually, $\Xi=\{1,\ldots,n\}\times\{1,\ldots,n\}$).

Every nonnegative square matrix $A=(a_{ij})_{i,j\in\calV}$ corresponds to a weighted digraph $\calG[A]=(\calV,\calE[A],A)$ whose nodes are indexed by $\calV$ and whose set of arcs is defined as $\calE[A]=\{(i,j)\in\calV\times\calV:a_{ij}>0\}$. The value $a_{ij}$ can be interpreted as the weight of arc $i\rightarrow j$. A sequence of arcs $i_0\rightarrow i_1\rightarrow\ldots\rightarrow i_{s-1}\rightarrow i_s$ constitute a \emph{walk} between nodes $i_0$ and $i_s$ in graph $\calG[A]$. The set of nodes $J\subseteq\calV$  is \emph{reachable} from node $i$ if $i\in J$ or a walk from $i$ to some element $j\in J$ exists; $J$ is called \emph{globally reachable} in the graph if it is reachable from every node $i\not\in J$.

A graph is strongly connected (strong) if every two nodes $i,j$ are mutually reachable. A graph that is not strong has several strongly connected (or simply strong) \emph{components}. A strong component is said to be non-trivial if it contains more than one node. A component is said to be a \emph{sink} component if no arc leaves it and a \emph{source} component if no arc enters it.
A strong component can be \emph{isolated}, when it has neither incoming nor outcoming arcs, and thus it is both a source and a sink.
Strong components of undirected graphs are always isolated.
\begin{figure}[h]
\begin{subfigure}{0.49\columnwidth}
\centering
\includegraphics[width=\textwidth]{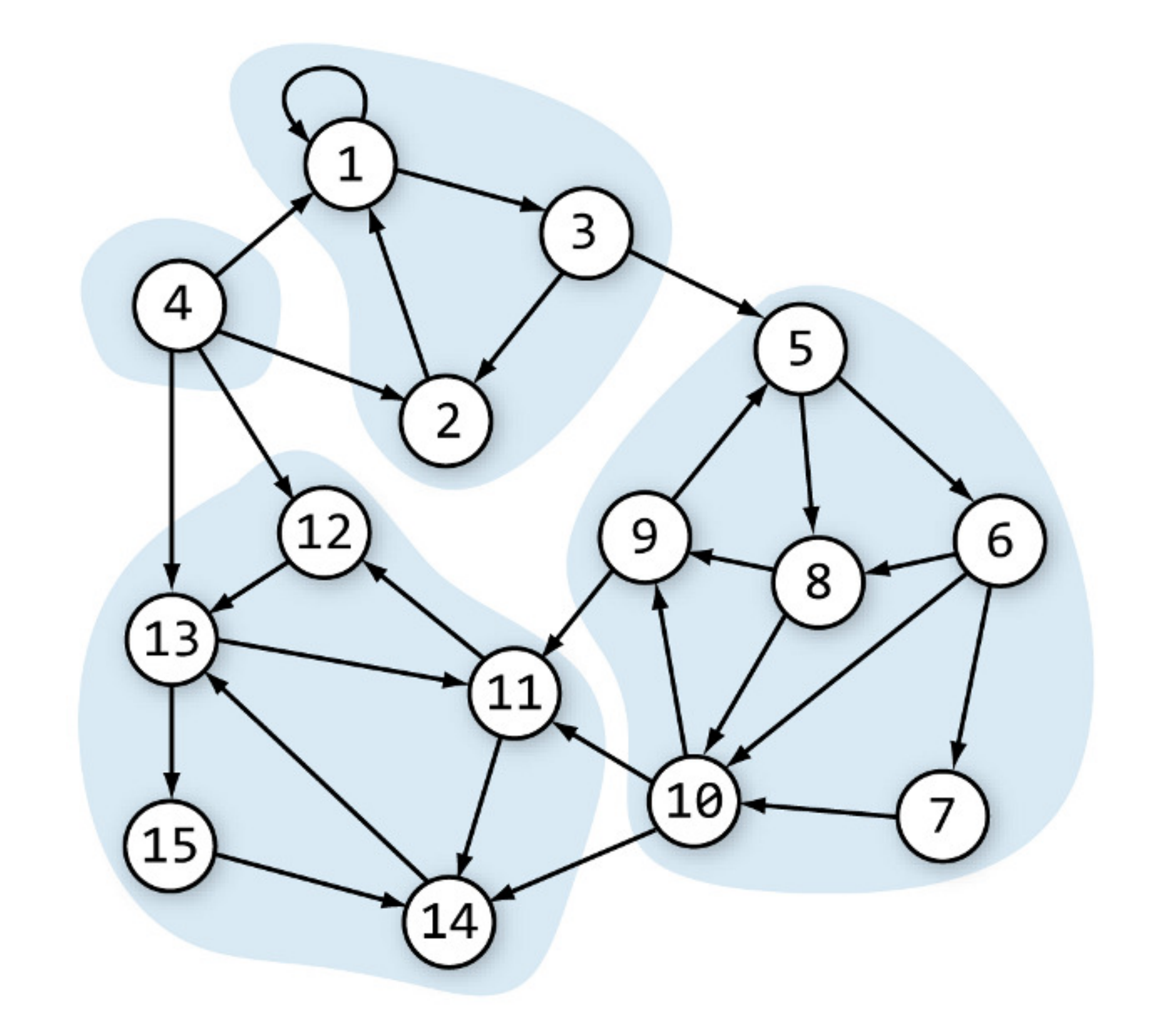}
\caption{}
\end{subfigure}
\begin{subfigure}{0.49\columnwidth}
\centering
\includegraphics[width=\textwidth]{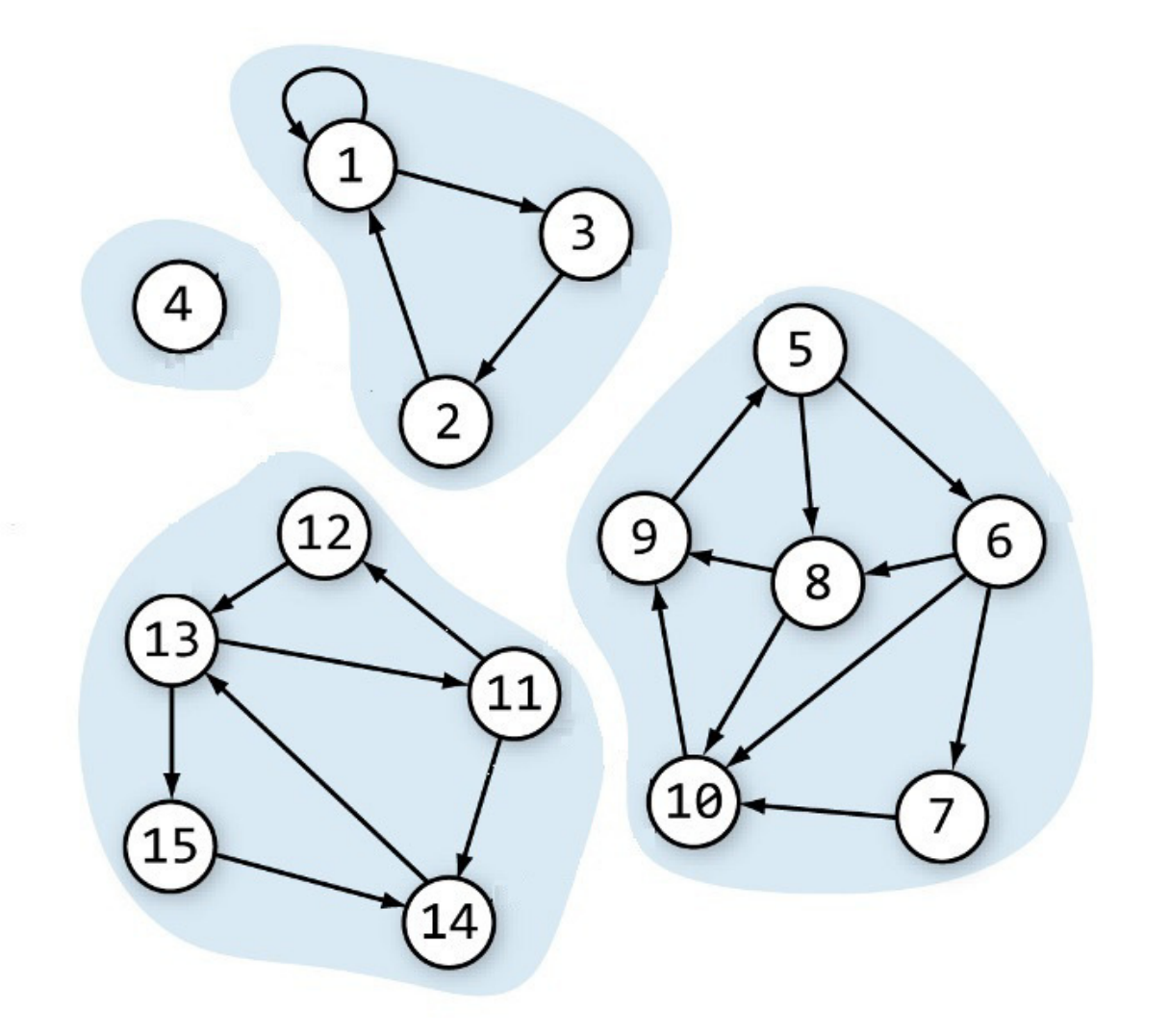}
\caption{}
\end{subfigure}
\caption{Strong components of a directed graph: (a) non-isolated; (b) isolated. In (a), $\{4\}$ is a (trivial) single source component, $\{11,\ldots,15\}$ is a single sink component.}
\label{fig.strong}
\end{figure}


\section{The Eisenberg-Noe financial network model}
\label{sec:oneperiod}
%
We start by considering  the ``static'' case introduced in the seminal work of Eisenberg and Noe \citep{EisNoe:01}.
In this setting,  $n$ nodes, representing financial entities (banks),  are connected via a complex
structure of  mutual liabilities.
The payment due from node $i$ to node $j$ is denoted by $\bar p_{ij} \geq 0$, and such liabilities are supposed to be due at the end of a fixed time period.
These interbank liabilities form the {\em liability matrix}
$\bar P	\in\Real{n\times n}$, such that $[\bar P]_{ij} = \bar p_{ij} $ for $i\neq j=1,\ldots,n$, and
$[\bar P]_{ii} = 0$ for $i=1,\ldots,n$.

Following the notation introduced in \cite[Section~5]{GlasYou:16}, we let $c\in\Realp{n}$ be the vector whose $i$th component $c_i\geq 0$ represents the total payments due to node $i$ from non-financial entities (i.e., from any other entity, different from the $n$ banks). Payments from banks to the external sector are instead modeled
%
by introducing a fictitious node  that represents the external sector and owes no liability to the other nodes (the corresponding row of $\bar P$ is zero).

The nominal cash in-flow and out-flow at a node $i$ are, respectively,
\[
\bar\phi_i\ap{in}\doteq  c_i + \sum\nolimits_{k\neq i} \bar p_{ki},\quad
\bar p_i \doteq  \bar\phi_i\ap{out}\doteq   \sum\nolimits_{k\neq i} \bar p_{ik}.
\]
In regular operations, the in-flow at each bank is no smaller than its out-flow (i.e.,  $\bar\phi_i\ap{in} \geq \bar\phi_i\ap{out}$),
each bank remains solvable and is able to pay its liabilities in full.
%
A critical
 situation occurs instead when (due to, e.g., a drop in the external liquidity in-flow $c_i$) some bank $i$
 has not enough incoming  liquidity to fully  pay its liabilities.
 In this situation, the actual payments to other banks  have to be remodulated to lesser values than  their nominal values $\bar p_{ij}$.
 The {\em clearing payments} are a set of mutual payments which settle the mutual claims in case of defaults, by enforcing a set of rules~\citep{EisNoe:01,csoka2018decentralized}, which are: {\em (i)} payments cannot exceed the corresponding liabilities,
 {\em (ii)} {\em limited liability}, i.e., the balance at each node 
 cannot be negative,
 {\em (iii)} {\em absolute priority} (i.e., each node either pays its liabilities in full, or it pays out all its balance).

We let $p_{ij} \in [0,\bar p_{ij}]$, $i\neq j=1,\ldots,n$, denote the actual inter-bank payments executed at the end of the period, which we shall collect in matrix $P\in\Real{n,n}$.
At each node $i$ we write a flow balance equation, involving the actual cash in-flow and  out-flow, defined respectively as
\begin{gather}
\phi_i\ap{in} \doteq c_i + \sum\nolimits_{k\neq i} p_{ki} ,\label{eq:in-flow}\\
\phi_i\ap{out} \doteq  p_i \doteq  \sum\nolimits_{k\neq i} p_{ik}.\label{eq:out-flow}
\end{gather}
 The  cash balance represents the net worth $w_i$ of the $i$th bank, which is defined as
 \beq
 w_i \doteq \phi_i\ap{in} - \phi_i\ap{out} =  c_i + \sum\nolimits_{k\neq i} p_{ki} -  \sum\nolimits_{k\neq i} p_{ik}.
\label{eq:networth}
\eeq
The limited liability rule {\em (ii)} requires that $w_i\geq 0$, $\forall i$.

In vector notation, the vectors of actual and nominal in/out-flows and the vector of net worths are
\begin{gather}
 \phi\ap{in} = c + P\tran \one,\quad
 \bar \phi\ap{in} = c + \bar P\tran \one \label{eq:in-flow-vec}\\
 \phi\ap{out} = p =  P \one,\quad
  \bar\phi\ap{out} =  \bar p = \bar P \one\label{eq:out-flow-vec}\\
  w = \phi\ap{in} - \phi\ap{out} = (c+P\tran \one )- P\one,\label{eq:ntworth}
   \end{gather}
where $\one$ denotes a  vector of ones of suitable dimension.

The above mentioned conditions  {\em (i), (ii)}  on the payments are written in compact vector form
as $0\leq P \leq \bar P$ and $P\one \leq c+P\tran \one$, that is
the payment matrix $P$ is restricted to belong to the following convex polytope
\beq\label{eq.polytope-matr}
\begin{aligned}
\calP(c,\bar P)& \doteq \left\{P\in\Real{n\times n}:\, 0 \leq  P\leq \bar P,\right.\\
&\left.P\one \leq c+P\tran \one,\, P_{ii}=0,\;i = 1,\ldots,n
\right\}.
\end{aligned}
\eeq
A payment matrix  $P\in\calP(c,\bar P)$ is a {\em clearing matrix}, or matrix of clearing payments, if it complies with the
 \emph{absolute priority} of debt claims rule {\em (iii)}, that is,
 \beq\label{eq:maximality}
 P\one = \min( \bar P\one, c+P\tran \one).
\eeq
It can be shown~\citep{csoka2018decentralized,Journal2021} that a clearing matrix
can be found by solving an optimization problem of the form
\beq \label{eq:clearingoptP}
\begin{aligned}
\min_{P} & \quad f(P)\\
\mbox{subject to:} &\quad P\in \calP(c,\bar p)
\end{aligned}
\eeq
where $f$ is a decreasing function of the matrix argument $P$ on $[0,\bar P]$, i.e., a function such that $\bar P\geq P^{(2)} > P^{(1)}\geq 0$, $P^{(2)} \neq P^{(1)}$, implies
$f(P^{(2)}) < f( P^{(1)})$. It can be shown that for any choice of $f$ the solution to~\eqref{eq:clearingoptP} is automatically a clearing matrix, that is,~\eqref{eq:maximality} holds.
Possible choices for  $f$  in~\eqref{eq:clearingoptP} are for instance
$f(P) = \|\bar \phi\ap{in} -\phi\ap{in}\|_1$ and
$f(P) = \|\bar \phi\ap{in} -\phi\ap{in}\|_2^2$, where
 $\phi\ap{in}(P)  = c +P\tran \one$. The optimal solution of \eqref{eq:clearingoptP}, however,
 may be non unique in general~\citep{Journal2021}.

 \subsection{The pro-rata rule}\label{subsec.prorata}
 In practice, payments under default are subject to additional prevailing regulations.
 A common one is the so called
 {\em proportionality} (or, pro-rata) rule, according to which payments are made in proportion to the original outstanding claims.
 Denoting by
 \beq
 a_{ij}\doteq \left\{\ba{ll}\dss \frac{\bar p_{ij}}{\bar p_i}  & \mbox{if }  \bar p_i> 0 \\
 1 & \mbox{if }  \bar p_i=0\,\mbox{ and } i=j \\
 0 & \mbox{otherwise}
 \ea\right.
 \label{eq:prorataA}
 \eeq
  the relative proportion of payment due nominally by node $i$ to node $j$,  we  form the {\em relative liability} matrix $A = [ a_{ij}]$.
     By definition, $A$ is row-stochastic, that is  $A\one = \one$.
The pro-rata rule imposes  the relations
 \beq
  p_{ij} = a_{ij}  p_i ,\quad \forall i,j,
  \label{eq_prorata}
 \eeq
 where $p_i$ is the out-flow defined in~\eqref{eq:out-flow}.  
 In matrix notation, the pro-rata rule corresponds to a linear equality constraint on the entries of $P$, that is
 $P = \mbox{diag}(P\one) A$. Under pro-rata rule, the problem of clearing payments can be rewritten in terms of  the total out-payments vector $p = P\one$, which is said to be feasible if it belongs to
 \beq\label{eq.polytope-vect}
 \calP\ped{pr}(c,\bar p) \doteq \{p\in\Real{n}:\, 0 \leq  p\leq \bar p ,\; p \leq c+ A\tran p\},
 \eeq
 where $\bar p\doteq \bar P\one$. Among the feasible payment vectors $p\in\calP\ped{pr}(c,\bar p) $, a vector of clearing payments, or simply {\em clearing vector} is a vector $p\in\calP\ped{pr}(c,\bar p) $ such that
 \[
 p = \min( \bar p, c+A\tran p).\tag{\ref{eq:maximality}a}\label{eq:maximalitya}
\]
A clearing vector  $p^*$ can be found~\citep{GlasYou:16} by solving an optimization problem of the form
\beq\label{eq:clearingopt}
\begin{aligned}
\min_{p} & \; f(p) \\
\mbox{subject to:} &\quad p\in \calP\ped{pr}(c,\bar p)
\end{aligned}
\eeq
where $f:[0,\bar p]\to\Real{}$ is any decreasing function,  that is, a function such that
 $p^{(1)},p^{(2)}\in[0,\bar p]$ and $p^{(1)}\leq p^{(2)}$ imply  $f(p^{(1)})\geq f(p^{(2)})$, and the latter inequality is strict unless $p^{(1)}=p^{(2)}$.
 Possible choices for  $f$ are for instance $f(p ) = \|\bar \phi\ap{in} -\phi\ap{in}(p) \|_2^2$,  and $f(p ) =\sum_{i=1}^n (\bar \phi\ap{in}_i -\phi\ap{in}_i(p) )$, where  $\phi\ap{in}(p)  = c +A\tran p$.
The following proposition holds.

 \begin{lemma}\label{lem.unique}\cite[Lemma~1]{Journal2021}
 The solution $p^*=p^*[A,c,\bar p]$ to~\eqref{eq:clearingopt} exists, is unique and does not depend on the choice of
 $f$, provided that $f$ is decreasing. Additionally,
 \begin{enumerate}[(a)]
 \item $p^*$ satisfies~\eqref{eq:maximalitya} (being thus a clearing vector);
 \item $p^*$ dominates any other admissible payment vector
 \[
 p^*[A,c,\bar p])\geq p\;\;\forall p\in \calP\ped{pr}(c,\bar p);
 \]
 \item each strongly connected component being a sink (without out-coming arcs) of graph $\calG[A]$ contains at least one node $i$ such that $p_i^*=\bar p_i$;
\item $p^*$ is the unique solution of ~\eqref{eq:maximalitya} enjoying the property from statement (c);
 \end{enumerate}
 \end{lemma}

 Lemma~\ref{lem.unique}, in fact, is valid for \emph{every} stochastic matrix $A\in\Real{\calV\times\calV}$, because its proof (available in~\cite{Journal2021}) does not rely on~\eqref{eq:prorataA}.

 \section{Dynamic financial networks}\label{sec:dynamics}

 A key observation is that the  default and clearing model discussed in the previous section, which coincides with the mainstream one studied in the literature~\citep{GlasYou:16}
 is an {\em instantaneous} one. By instantaneous we  mean that the described process assumes that at one point in time (say, at the end of a day), all liabilities are claimed and due simultaneously, and that the entire network of banks becomes aware of the claims and possible defaults and instantaneously agrees on the clearing payments.
 On the one hand such an instantaneous model may be quite unrealistic, and on the other hand
  the implied default mechanism is such that all financial operations of  defaulted nodes are instantaneously frozen, which possibly induces propagation of the default to other neighboring nodes, in an avalanche fashion, see, e.g.~\cite{Massai:2021}.

One motivation for the dynamic model we propose in this paper is that one may expect that if financial operations are allowed for a given number of time periods after the initial theoretical defaults, some nodes may actually {\em recover} and eventually manage to fulfill their obligations. The overall system-level advantage of such strategy is that the catastrophic effects of avalanche defaults are possibly mitigated, as shown by examples in Section~\ref{sec.exam}.

In our dynamic multi-period model described below, if a theoretical default condition (we shall call this a {\em pseudo-default}) happens at some time $t < T$, where $T$ is the final time, we do not freeze operations. Instead, we carry over the residual liabilities for the next period and let the nodes continue their mutual payments operations, and so on until the final time $T$.
The key elements of this model are the following:
\begin{itemize}
 \item $t=0,1,\ldots,T$, denote discrete time instants  delimiting periods of fixed length (e.g, one day, one month, etc.);
 \item $T\geq 0$ denotes the final horizon;
 \item $c(t) \in\Real{n}\geq 0$ represents the cash in-flow at the nodes at the beginning of period $t$;
 \item  matrix $\bar P(t)=(\bar p_{ij}(t))\in\Real{n,n}$ describes the liabilities (i.e., the mutual payment obligations) among the nodes at period $t$, i.e., $\bar p_{ij}(t)$ is the nominal amount due from $i$ to $j$ at the end of period $t$.
 $ \bar P \doteq \bar P(0)$ denotes the initial liabilities at $t=0$;
 \item matrix $P(t)=(p_{ij}(t))\in\Real{n,n}$ contains the actual payments from $i$ to $j$ performed at the end of period $t$;
 \item the vectors of actual and nominal in-flows and out-flows $\phi\ap{in} (t),\phi\ap{out} (t),\bar\phi\ap{in} (t),\bar\phi\ap{out} (t)$ at period $t=0,\ldots,T-1$, are defined similarly to~\eqref{eq:in-flow-vec} and~\eqref{eq:out-flow-vec};
 \item the net worth $w_i(t)$ of node $i$ at the beginning of period $t$ evolves in accordance with
 \beq
 w_i(t+1) = w_i(t) + \phi_i\ap{in} (t)  - \phi_i\ap{out}(t)
 \label{eq_wdynamics}
 \eeq
 or, in the equivalent vector form
 \beq
 w(t+1) = w(t) + c (t) + P(t)\tran \one - P(t) \one. 
 \label{eq_wdynamics2_free}
\eeq
 \end{itemize}
Similar to the single-period case discussed in Section~\ref{sec:oneperiod}, the limited liability condition requires that
$w(t)\geq 0$ at all $t$.
It may therefore happen that
 a payment $p_{ij}(t)$ has to be lower than the corresponding liability $\bar p_{ij}(t)$ in order to guarantee
 $w_i(t) \geq 0$.
When this happens at some $t< T$, instead of declaring default and freezing the financial system, we allow operations to continue up to the final time $T$, updating the
due payments according to the  equation
  \beq
\bar p_{ij}(t+1) =  \alpha \left(\bar p_{ij} (t) - p_{ij} (t)  \right) ,
  \label{eq_pdynamics}
 \eeq
 where $\alpha \geq 1$ is the interest rate applied on past due payments.
The previous relation  can be written as
 \beq
\bar P(t+1) = \alpha \left(\bar P (t) - P (t)  \right), \quad t\in \calT,
  \label{eq_pdynamics2}
 \eeq
 where $\calT \doteq \{0,\ldots,T-1\}$.   The meaning of equation~\eqref{eq_pdynamics2} is that if  a due payment at $t$ is not paid in full, then the residual debt is added to the nominal liability for the next period,
 possibly increased by an interest factor $\alpha \geq 1$. This mechanism allows for a node which is technically in default at a time $t$ to continue operations and (possibly) repay its dues in subsequent periods.
 Notice that time-varying $\bar P(t)$ depends on the \emph{actual} payment matrices $P(0),\ldots,P(t-1)$.
 The final nominal matrix $\bar P(T)$ contains the residual debts at the end of the final period.
 The recursions~\eqref{eq_wdynamics2_free} and~\eqref{eq_pdynamics2} are
  initialized with
 \beq
 w(0)  = 0,\quad \bar P(0 ) =  \bar P,
 \label{eq:initaliz}
 \eeq
where  $\bar P$ is the initial liability matrix.

 %
 %

 Vectors of external payments $c(t)$  are considered as given inputs, while
actual payments matrices $P(t)$ are to be determined, being subject to the constraints
\begin{gather}
P(t) \geq 0, \quad P(t) \leq \bar P(t), \quad t\in \calT \label{eq.cond-p-1}\\
 P(t) \one \leq   w(t) + c (t) + P(t)\tran \one , \quad  t\in \calT,\label{eq.cond-p-2}
\end{gather}
where \eqref{eq.cond-p-1} represents the requirement that actual payments never exceed the nominal liabilities, and
\eqref{eq.cond-p-2} represents the requirement that $w(t+1)$, as given in \eqref{eq_wdynamics2_free}, remains nonnegative at all $t$.
%
Conditions  \eqref{eq.cond-p-1}, \eqref{eq.cond-p-2}  can be made explicit by eliminating the variables $w(t)$ and $\bar P(t)$, which
by using \eqref{eq_wdynamics2_free}--\eqref{eq:initaliz} can be expressed as
\begin{gather}
 \bar P(t) = \alpha^t \bar P(0) -  \sum\nolimits_{k=0}^{t-1} \alpha^{t-k}P(k), \label{eq.bar-p-evolves}\\
 w(t) = C(t-1) + \sum\nolimits_{k=0}^{t-1} \left(P\tran(k) - P(k) \right)\one , \label{eq.w-evolves}\\
C(t) \doteq \sum\nolimits_{k=0}^{t} c(k),\quad t=0,\ldots,T .  \label{eq:cumulativeinflow}
\end{gather}
Conditions~\eqref{eq.cond-p-1}, \eqref{eq.cond-p-2} can thus be rewritten as
\begin{gather}
P(t)\geq 0,\label{eq.cond-p-1+}\\
\sum\nolimits_{k=0}^{t} \alpha^{t-k}P(k) \leq \alpha^t \bar P\label{eq.cond-p-1a}\\
C(t) + \sum\nolimits_{k=0}^{t} \left(P(k) \tran- P(k) \right)\one \geq 0\label{eq.cond-p-2a}\\
\forall t\in \calT \notag.
\end{gather}
For brevity, we denote
\[
[P]\doteq (P(0),\ldots,P(T-1)),\; [c]\doteq (c(0),\ldots,c(T-1)).
\]

\begin{defn}
We call a sequence of payment matrices $[P]$ \emph{admissible} if conditions
\eqref{eq.cond-p-1+}--\eqref{eq.cond-p-2a}  hold. Let
\[
\mathcal{P}([c],\bar P)\doteq\{[P]:\text{\eqref{eq.cond-p-1+}--\eqref{eq.cond-p-2a} hold}\}
\]
stand for the polyhedral set of all admissible matrix sequences $[P]$ that correspond to the given sequence of vectors $[c]$
and initial liability matrix $\bar P$.
\end{defn}

The system-level cost  that we consider is  the cumulative sum of deviations
of the actual in-flows at nodes from the nominal ones, that is
\beq\label{eq.cost}
L([P])  \doteq  \sum_{t=0}^{T-1}  \sum_{i=1}^n  (\bar \phi\ap{in}_i(t) -  \phi\ap{in}_i(t) ).
\eeq
From the definition~\eqref{eq:in-flow-vec} of in-flow vectors and from \eqref{eq.bar-p-evolves} we obtain that
\bea
L([P])  &=&  \sum_{t=0}^{T-1} \one\tran  ( \bar \phi\ap{in}(t) -  \phi\ap{in}(t)) =
\sum_{t=0}^{T-1} \one\tran  (\bar P(t)-P(t))\one \nonumber \\
&=& \sum_{t=0}^{T-1} \one\tran (  \alpha^t \bar P  -\sum_{k=0}^t \alpha^{t-k} P(k) ) \one \nonumber \\
&=& a_0 \one \tran \bar P  \one - \sum_{t=0}^{T-1} a_t \one\tran P(t)  \one,\nonumber
\eea
where  the constants $a_0>a_1>\ldots>a_{T-1}$ are defined as
\beq
a_t\doteq \sum_{j=0}^{T-t-1} \alpha^j = \left\{ \ba{ll}
 \frac{\alpha^{T-t}-1}{\alpha -1}, & \mbox{if } \alpha > 1\\
T-t, & \mbox{if } \alpha = 1.
\ea \right.
\label{eq:ai}
\eeq
The optimal payment matrices are thus obtained as a solution to the following optimization problem
\beq\label{eq:prob01a}
\max_{[P]} \sum_{t=0}^{T-1} a_t \one\tran P(t)  \one\quad \mbox{s.t.:} \quad [P]\in\calP([c],\bar P),
\eeq
which is equivalent to minimization of the overall ``system loss'' $L([P])$ over the set of all admissible payment matrices.

Observe that, from a numerical point of view,
finding an optimal sequence of payment matrices amounts to solving the linear programming (LP) problem \eqref{eq:prob01a}.
Notice also that in the case $T=1$ the set $\mathcal{P}([c])$ reduces to the polytope of matrices~\eqref{eq.polytope-matr}, and the optimization problem~\eqref{eq:prob01a}
is a special case of~\eqref{eq:clearingoptP}, where $f(P)=- \one\tran P(0)  \one$. 

We next establish a fundamental property of the payment matrices resulting from \eqref{eq:prob01a}.

\subsection{The absolute priority rule}
Recall that in the static (single period) case the optimal payment matrix automatically satisfies the absolute priority rule~\eqref{eq:maximality}. A natural question arises whether a counterpart of this rule
can be proved for the dynamical model in question: is it true that a bank failing to meet the nominal obligation has to
nevertheless pay the maximal possible amount? Mathematically, this means that for all $t=0,\ldots,T-1$
the following implication holds:
\beq\label{eq.default-imp}
\phi_i\ap{out}(t)<\bar\phi_i\ap{out}(t)\Longrightarrow \phi_i\ap{out}(t)=\phi_i\ap{in}(t)+w(t).
\eeq
The affirmative answer is given by the following theorem.

\begin{theorem}\label{thm.1}
Suppose that $[P]=(P(t))_{t=0}^{T-1}$ is an {optimal} solution of~\eqref{eq:prob01a}, and let $(\bar P(t))_{t=0}^T$ be the corresponding sequence of nominal liability matrices, defined in accordance to~\eqref{eq_pdynamics2}.  For a given bank $i$, let $t_*=t_*(i)$ be the {first} instant when $i$ pays its debt to the other banks
\[
p_{ij}(t_*)=\bar p_{ij}(t_*)\quad\forall j\ne i
\]
(if such an instant fails to exist, we formally define $t_*=T$).
Then, either $t_*=0$ (the debt is paid immediately) or
\beq\label{eq.default}
\phi_i\ap{out}(t)=\phi_i\ap{in}(t)+w_i(t)\;\;\forall t = 0,\ldots,(t_*-1).
\eeq
In particular, the implication~\eqref{eq.default-imp} holds for any optimal sequence of payments matrices $[P]$.
Furthermore, for each $t\geq 1$ the graph $\calG[P(t)]$ contains no directed cycles.
\end{theorem}

A proof of Theorem~\ref{thm.1} is provided in Appendix~\ref{sec.proofs_th1}.

\begin{remark}\label{rem.earliest}\rm
Implication~\eqref{eq.default-imp} implies that each bank pays its nominal liability
at the earliest period $t$ when such a payment is possible: $w_i(t)+\phi_i\ap{in}(t)\geq \bar\phi_i\ap{out}(t)$.
The requirement of minimal system loss prevents unnecessary deferral of payments and pushes the banks towards paying the claims as early as possible.
Since the payment matrices resulting from the solution of \eqref{eq:prob01a}
satisfy the rules {\em (i), (ii), (iii)} from Section~\ref{sec:oneperiod}, they are guaranteed to be proper clearing matrices at each stage.
\erem
\end{remark}

\subsection{A sub-optimal sequential approach}\label{subsec.sequential}
\label{sec:sequential:full}
Looking at the objective function in problem  \eqref{eq:prob01a}, we observe that this function is linear and separable in the
$P(t)$ variables, $t=0,\ldots,T-1$. 
Also, looking at the constraints of \eqref{eq:prob01a}, given by  \eqref{eq.cond-p-1+}--\eqref{eq.cond-p-2a},
we see that at each $t=0,\ldots,T-1$ the variable $P(t)$ is constrained as
\[
\begin{gathered}
0 \leq P(t) \leq \bar P(t), \\
w(t) +  c(t) +  \left(P\tran(t) - P(t) \right)\one  \geq 0,
\end{gathered}
\]
where
\bea
\bar P(t) & \doteq & \alpha^t \bar P -  \sum_{k=0}^{t-1} \alpha^{t-k}P(k) \label{eq:Qeq}\\
w(t) &\doteq & C(t-1) +  \sum_{k=0}^{t-1} \left(P\tran(k) - P(k) \right)\one ,\label{eq:qeq}
\eea
and $\bar P(t) $, $w(t)$ depend only on the variables $P(0),\ldots,P(t-1)$
 and external payments $c(0),\ldots,c(t-1) $ at periods preceeding $t$.
 This suggests the following recursive relaxation of problem \eqref{eq:prob01a} where, at each $t=0,\ldots,T-1$
 we solve a problem in the $P(t)$ variable only
 \bea
\tilde P^*(t) = \arg\max_{P(t)} && \one\tran P(t)  \one  \label{eq_recursivep} \\
\mbox{s.t.:} && w^*(t) +  c(t) +  \left(P\tran(t) - P(t) \right)\one  \geq 0 , \nonumber \\
&&  0 \leq P(t) \leq \bar P^*(t-1) ,\nonumber
\eea
where
 $\bar P^*(t)$, $w^*(t) $ are given by \eqref{eq:Qeq}, \eqref{eq:qeq} evaluated at the previous optimal values
$\tilde P^*(0) , \ldots, \tilde P^*(t-1) $, and initialized so that
$w^*(0) \doteq  0 $, $\bar P^*(0) \doteq  \bar P$.

It is clear by construction that  any optimal sequence of solutions $\tilde P^*(0), \ldots,  \tilde P^*(T-1)$ of
\eqref{eq_recursivep} is feasible for problem \eqref{eq:prob01a}. However, this
 ``greedy'' sequential solution is in general not  optimal for problem   \eqref{eq:prob01a},
 as highlighted by the following example.

\begin{figure}[h]
\centering
\includegraphics[width=.4\textwidth]{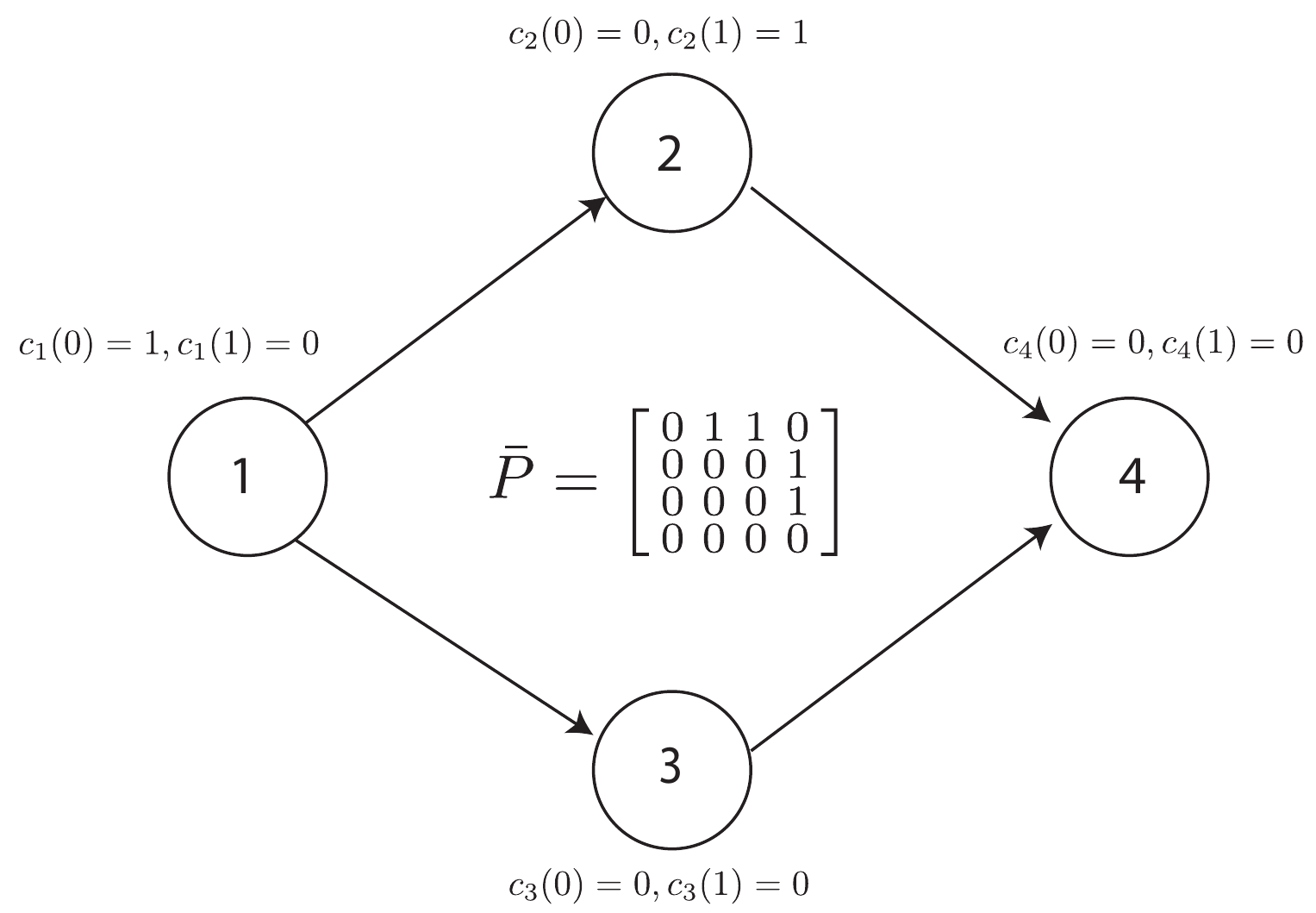}
\caption{A four-node liability network.}
\label{fig:counterex}
\end{figure}

\textbf{Example~1.} Consider a group of four banks with initial  liability matrix $\bar P$ and liability graph shown in Fig.~\ref{fig:counterex}.

We assume that $\alpha=1$ and consider a time horizon $T=2$, with external payments
 $c(0)=(1,0,0,0)^{\top}$,
$c(1)=(0,1,0,0)^{\top}$.
The unique optimal strategy in~\eqref{eq:prob01a} can be easily found: at stage $0$, node $1$ pays its maximum possible to node $3$, i.e., $p_{13}(0)=\bar p_{13}=1$, and node $3$ transfers it
 to node $4$: $p_{34}(0)=\bar p_{34}=1$.
 Node $2$ receives and pays nothing at period $t=0$, while at  $t=1$  node $2$ receives an external payment and hence  pays
 its liability to node $4$: $p_{24}(1)=\bar p_{24}=1$.  This optimal strategy leads to the  optimal loss $L=3$,
 and at the end of the time horizon only
 node $1$ is in default (owing $1$ to node $2$).
 %

If we consider the sequential approach instead, we see that the objective function~\eqref{eq_recursivep} at $t=0$
is $(p_{12}(0)+p_{13}(0)) +  p_{24}(0) + p_{34}(0)$,
%
hence it is insensitive to  how node $1$ divides its asset $c_1(0)=1$ between nodes $2$  and $3$.
An optimal solution to~\eqref{eq_recursivep} at $t=0$ is for instance
$\tilde p_{12}^*(0)=\bar p_{12}=1$,
$\tilde p_{24}^*(0)=\bar p_{24}=1$.
With this solution in place, problem \eqref{eq_recursivep} at $t=1$
leads to a network in which no further payments can be made (i.e., $\tilde p_{ij}^*(1) = 0$ for all $i,j$), and the loss function under this sub-optimal solution is $L=4$, with two defaulted nodes at the end of the horizon:
 node $1$, which still owes $1$ to node $3$, and node $3$, which still owes $1$ to node $4$.

The point here is that the correct choice at $t=0$ cannot be made in general unless one knows the {\em future} external payments at all nodes and at all $t>0$. The sequential solution hence remains sub-optimal, since it does not exploit this information (it only uses, at each $t$, the observed external payments $c(k)$, $k=0,\ldots,t$, up to that $t$).
On the one hand, this fact highlights that the solution to the ``full'' problem~\eqref{eq:prob01a} is in general superior
in terms of optimal loss to the solution of the sequential problem. On the other hand, however, it also underlines that the whole stream of future external payments must be known at $t=0$ in order to being able to solve~\eqref{eq:prob01a}. If, at each $t$, one has total uncertainty about the future payments $c(\tau)$, $\tau > t$, then the full approach is not viable
while the sequential approach still is.
\est

\subsection{Dynamic networks with pro-rated payments}

The pro-rata rule discussed in Subsection~\ref{subsec.prorata} can be introduced also in the dynamic network setting.
Here, we let the pro-rata matrix be fixed according to the {\em initial} liabilities, that is the $A$ matrix is given by \eqref{eq:prorataA} with $\bar P = \bar P(0)$.
Then, the pro-rata rule is nothing but a linear equality constraint on the payment matrices, that is
\beq
P(t) = \mbox{diag} (P(t)\one) A, \quad t=0,\ldots,T-1.
\label{eq:proratac}
\eeq
  In view of the definition of $A$, one has $\bar P(0)=\mbox{diag} (\bar P(0)\one) A$.
  Using induction on $t$ and equation~\eqref{eq.bar-p-evolves},  it can be easily shown that~\eqref{eq:proratac} entails the equations
  \[
  \bar P(t)=\mbox{diag} (\bar P(t)\one) A,\quad t=0,\ldots,T.
  \]
Hence, payment matrices $P(t)$ and $\bar P(t)$ are uniquely determined by the actual and nominal payment vectors
 \beq\label{eq.p-vs-phi}
 p(t) \doteq P(t) \one=\phi\ap{out}(t),\quad \bar p(t)\doteq\bar P(t)\one=\bar\phi\ap{out}(t).
\eeq
Also, it holds that $\phi\ap{in}=P\tran(t)\one = A\tran p(t)$.
Conditions~\eqref{eq.cond-p-1}, \eqref{eq.cond-p-2} can be now rewritten as
\begin{gather}
p(t)\geq 0,\label{eq.cond-p-1++}\\
\sum\nolimits_{k=0}^{t} \alpha^{t-k}p(k) \leq \alpha^t \bar p\label{eq.cond-p-1b}\\
C(t) + \sum\nolimits_{k=0}^{t} \left(A\tran p(k)- p(k) \right) \geq 0\label{eq.cond-p-2b}\\
\forall t\in\calT \notag.
\end{gather}

\begin{defn}
We call a sequence of payment vectors $[p]\doteq(p(0),\ldots,p(T-1))$ \emph{admissible} (under the pro-rata requirement) if conditions~\eqref{eq.cond-p-1++}--\eqref{eq.cond-p-2b}. Let
\[
\calP\ped{pr}([c],\bar p)\doteq\{[p]=(p(0),\ldots,p(T-1)):\text{\eqref{eq.cond-p-1++}--\eqref{eq.cond-p-2b} hold}\}
\]
stand for the convex polytope of all admissible sequences.
\end{defn}
Optimization problem~\eqref{eq:prob01a} can be now rewritten as
\beq
\max_{[p]} \sum_{k=0}^{T-1} a_k \one\tran p(k) \quad \mbox{s.t.:} \quad [p]\in\calP\ped{pr}([c],\bar p)\label{eq:prob01b}.
\eeq
This is again an LP problem, which may be solved numerically with great efficiency.
The pro-rata rule drastically reduces the number of unknown variables (each zero-diagonal payment $n\times n$ matrix reduces to $n$-dimensional vector). Furthermore,
unlike the original problem~\eqref{eq:prob01a}, the optimization problem~\eqref{eq:prob01b} admits a \emph{unique} maximizer $[p^*]$.
Also,  the solution abides by the absolute priority rule~\eqref{eq.default-imp}.
These properties are summarized in the following theorem.
\begin{theorem}\label{thm.prorata}
For each sequence $[c]$, the optimization problem~\eqref{eq:prob01b} has a {unique} solution $[p^*]$.
Furthermore, at each period $t=0,\ldots,T-1$, the optimal vector $p^*(t)$ is the unique solution of the LP:
\bea
 p^*(t)= \arg \max_{p}\;   1\tran p \label{eq:clearingopt-b} \rule{5cm}{0cm}\\
  \mbox{s.t.:}\;  0\leq p\leq\bar p^*(t),\, p\leq c(t)+w^*(t)+A\tran p,\label{eq:clearingopt-c} \rule{.6cm}{0cm}
\eea
where $ w^*(0)\doteq 0$, $ \bar p^*(0)\doteq \bar p$, and, for $t=1,\ldots,T-1$,
 \bea
 \bar p^*(t) & \doteq & \alpha^t \bar p -  \sum_{k=0}^{t-1} \alpha^{t-k} p^*(k) \label{eq:Qeq_pr}\\
w^*(t) &\doteq & C(t-1) +  \sum_{k=0}^{t-1} \left(A\tran  p^*(k) -  p^*(k) \right).\label{eq:qeq_pr}
\eea
In particular,  $p^*(t)\geq 0$ obeys the absolute priority rule
\beq
p^*(t) = \min(\bar p^*(t) , c(t) + w^*(t) +  A\tran p^*(t)).
\label{eq:dynfixpoint}
\eeq
\end{theorem}
The proof of Theorem~\ref{thm.prorata} is based on Lemma~\ref{lem.unique}, and it is
detailed in the Appendix~\ref{sec.proofs}.

\begin{remark}\rm
A few observations are in order regarding Theorem \ref{thm.prorata}.
First, we observe that the ``full'' multi-period problem \eqref{eq:prob01b} is
equivalent to the sequence of problems \eqref{eq:clearingopt-b}. Therefore, in  the pro-rata case
the sequential approach is \emph{optimal}, and not only sub-optimal, as it instead happened in the case with unrestricted payment matrices discussed in Section \ref{sec:sequential:full}.
Thus, the system-level objective in the full optimization problem~\eqref{eq:prob01b}
is minimized by finding regular clearing payments at each step $t$, whereby the liabilities among nodes are updated at each step by considering the  residual payments due to pseudo-defaults at the previous step.

Further,
we observe that, for each $t$,
problem~\eqref{eq:clearingopt-b}-\eqref{eq:clearingopt-c} has the same structure as
problem~\eqref{eq:clearingopt}, with $c=c(t)+w^*(t)$.
Hence, in view of the maximality of the vector $p^*(t)$, we have  that
the objective~\eqref{eq:clearingopt-b} can be replaced by any other \emph{increasing} function of $p$.
In view of Lemma~\ref{lem.unique}, the relations~\eqref{eq:clearingopt-b},\eqref{eq:clearingopt-c} can be rewritten as follows
\beq\label{eq.p*-dynamic}
p^*(t)=p^*[A,c(t)+w^*(t),\bar p^*(t)],
\eeq
which also entails~\eqref{eq:dynfixpoint} due to Lemma~\ref{lem.unique}, statement (a).
\erem
\end{remark}

 \section{Numerical illustration}\label{sec.exam}
 We consider a variation on the simplified network given in~\cite{GlasYou:16}.
  \begin{figure*}[t!h]
\centering
\includegraphics[width=.75\textwidth]{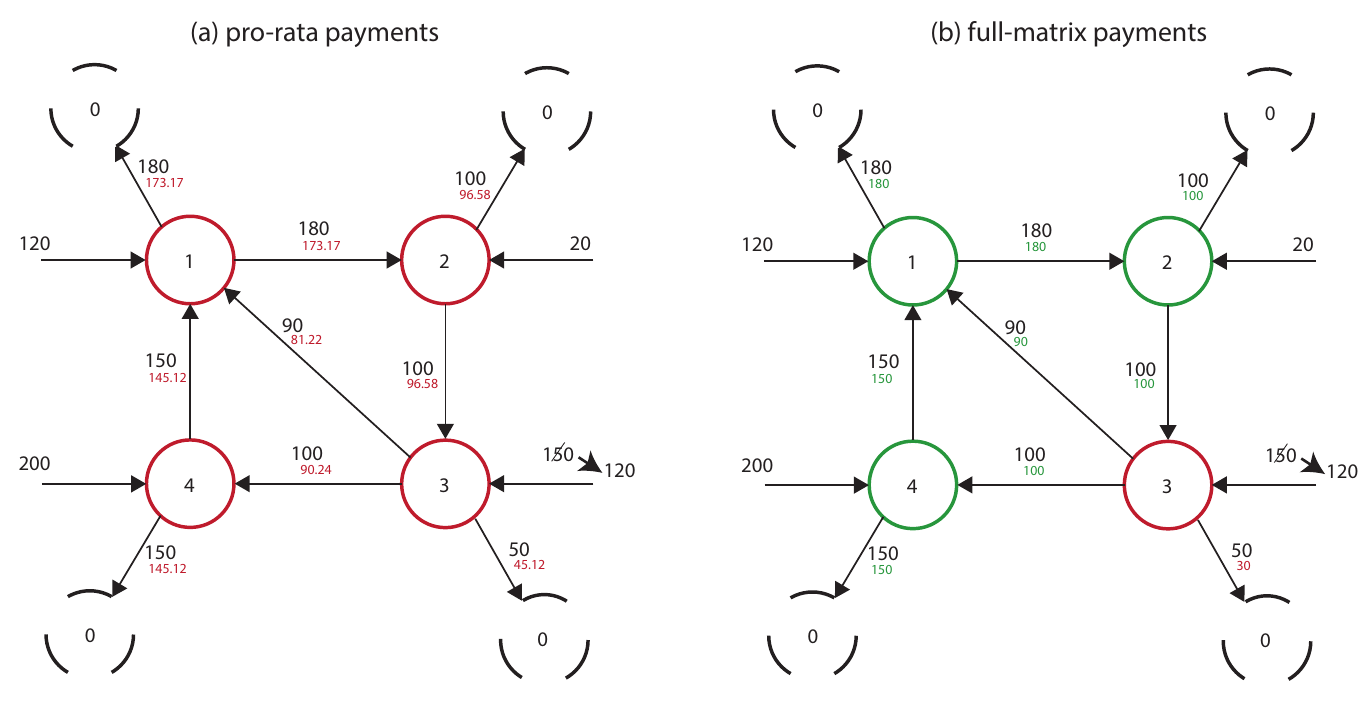}
\caption{Clearing payments in the example network. Left panel (a) shown the payments under pro-rata rule, Right panel (b)  shown the unrestricted clearing payments.}
\label{fig:clearings_ex1}
\end{figure*}
This network, displayed in Figure~\ref{fig:clearings_ex1},  contains $n=5$ nodes (including the fictitious sink node representing the external sector), with initial liability matrix
{\small
\[
\bar P =  \left[\begin{array}{ccccc} 0 & 180 & 0 & 0 & 180\\ 0 & 0 & 100 & 0 & 100\\ 90 & 0 & 0 & 100 & 50\\ 150 & 0 & 0 & 0 & 150\\ 0 & 0 & 0 & 0 & 0 \end{array}\right],
\]
}
where the last row refers to the  sink node.
We first discuss the static case, comparing pro-rata based results obtained by solving \eqref{eq:clearingopt} with those obtained
using an unrestricted payment matrix resulting from the solution of \eqref{eq:clearingoptP}.
Suppose there is a nominal  scenario where external cash flows are given as
\[
c = c\ped{nom} \doteq [120,\; 20,\; 150,\; 200,\; 0]\tran .
\]
It can be readily verified that in the nominal scenario all the nodes in the network remain solvent, and the clearing payments coincide with the nominal liabilities.
Consider next a situation in which  ``shock'' happens on the in-flow at node $3$, so that this in-flow reduces from $150$ to $120$, that is
  \[
c = c\ped{shock}\doteq  [120,\; 20,\; 120,\; 200,\; 0]\tran  .
\]
Under the pro-rata rule, the clearing payments, resulting from the solution of  \eqref{eq:clearingopt}, are shown in smaller font below the nominal liabilities in the left panel of Figure~\ref{fig:clearings_ex1}: all nodes in the network default in a cascade fashion due to initial default of node $3$. The total defaulted amount (the sum of all the unpaid liabilities) is in this case $53.66$.

Then, we dropped the pro-rata rule, and we computed the clearing payments
according to \eqref{eq:clearingoptP}. The results in this case are shown in the right panel of Figure~\ref{fig:clearings_ex1}:
only node $3$ defaults, while all other nodes manage to pay their full liabilities. Not only we reduced the  sum of all unpaid liabilities
 to $20$ (i.e., a $62.7\%$ decrease with respect to the pro-rata case), but we also obtained {\em isolation} of the contagion, since the default was
 confined to node $3$ and did not spread to other parts of the network.

We next considered the dynamic case. In both the pro-rata case and the full matrix case, the idea implied by the single-period (static) approach
is that in case of default the financial operations of a node are frozen, that is, defaulted nodes cannot operate even if there are cash in-flows that are foreseen in the immediate future. A classical situation arises when there is a liquidity crisis, i.e., due payments from the external sector are late and their lateness provokes defaults at some nodes, which freeze and may propagate further defaults over the network.
Suppose that the in-flows in the nominal vector $c\doteq [121,21,130,204,0]\tran$ do not arrive simultaneously at time $t=0$, due to delays, and the stream of in-flows is
\[
\begin{gathered}
c(0) = [60,10,120,0,0]\tran\\
c(1) = [60,8,0,200,0]\tran,\\
c(2) = [1,3, 10, 4, 0]\tran.
\end{gathered}
\]
\noindent
A static approach at time $t=0$, with unrestricted payment matrix, would result in
{\small
\[
  P = \begin{bmatrix} 0 & 180 & 0 & 0 & 70 \\ 0 & 0 & 100 & 0 & 90\\
 90 & 0 & 0 & 100 & 30 \\ 100 & 0 & 0 & 0 & 0\\ 0 & 0 & 0 & 0 & 0 \end{bmatrix},
\]
}

\noindent
 with all nodes in default and a total default loss of  $343.40$.
 If we allow operations to continue over an horizon $T=3$, according to the model described in Section~\ref{sec:dynamics}, assuming an interest rate $\alpha = 1.01$ (i.e., $1\%$ interest per period), and solving
 the multi-period problem \eqref{eq:prob01a} with full payment matrices, we obtain $P(0) = P$, and
 {\small
 \beas
 P(1) =
  \left[\begin{array}{ccccc} 0 & 0 & 0 & 0 & 110.5 \\
  0 & 0 & 0 & 0 & 8 \\ 0& 0 & 0 & 0
   & 0 \\ 50.5 & 0 & 0 & 0 & 149.5 \\ 0 & 0 & 0 & 0 & 0 \end{array}\right],\;
 P(2) = \left[\begin{array}{ccccc} 0 & 0 & 0 & 0 & 0.61 \\
  0 & 0 & 0 & 0 & 2.12 \\ 0& 0 & 0 & 0
   & 10 \\ 0 & 0 & 0 & 0 & 2.02 \\ 0 & 0 & 0 & 0 & 0 \end{array}\right].
 \eeas
 }

 \noindent
 After these three rounds of clearing payments, only node $3$ is in default, owing a residual $10.51$ to the external sector.
 Notice that, if we used the pro-rata rule, thus solving the multi-step problem \eqref{eq:prob01b}, we would obtain a different set of (pro-rata) clearing payments, leading to a final situation of default at all nodes, with  a total defaulted value of $21.07$.

 \section{Conclusions}
 In this paper we explored dynamic clearing mechanisms in financial networks, under both pro-rata payment rules and
 unrestricted matrix payments. In both cases, we proposed to compute the clearing payments as optimal solutions to suitable multi-stage linear optimization problems, namely problem \eqref{eq:prob01a} for the unrestricted case, and problem \eqref{eq:prob01b} in the pro-rata case.
 Theorem~\ref{thm.1} establishes some fundamental properties of the solution in the unrestricted case stating, in particular, that  payments are not unnecessarily delayed when they are feasible (absolute priority of debt claims), so that the solutions are indeed clearing matrices at each stage.  Unrestricted optimal payments, however, are possibly non-unique and need be computed in a centralized way, since  knowledge of the whole network structure is necessary.
  Theorem~\ref{thm.prorata} establishes instead  key properties of the optimal   pro-rated payments. The key fact is that the solution is in this case unique and, moreover, it can be computed by solving sequentially
  a series of LP problems~\eqref{eq:clearingopt-b}-\eqref{eq:clearingopt-c}.
  In turn, in under mild hypotheses (for instance, when graph $\calG[A]$ has a unique and globally reachable sink node),
  the LP solutions coincide with the solution of a  series of fixed-point equations of the form \eqref{eq:dynfixpoint}. These equations are uniquely solvable, see the discussion in~\cite{Journal2021}, and their unique solution can be found by means of a decentralized algorithm known as  the {\em fictitious default algorithm} of \cite{EisNoe:01}.  Hence, the optimal multi-stage payments in the pro-rata case can be obtained by decentralized iterations among neighboring nodes.
  Numerical investigations, see, e.g., \cite{Journal2021}, suggested that proportional payments may lead to severely suboptimal clearings, and may be a concurring cause of  cascaded defaults: removing the pro-rata rule, both in the static and in the dynamic case, generally improves the high-level objective of reducing the systemic effects of defaults.

\bibliographystyle{agsm}
\begin{small}
\bibliography{biblio}

@inproceedings{CDC2021,
    author={Calafiore, G.C. and Fracastoro, G. and Proskurnikov, A.V.},
    title={On Optimal Clearing Payments in Financial Networks (submitted)},
    year={2021},
    booktitle={IEEE Conf. Decision and Control},    
}

@article{Journal2021,
    author={Calafiore, G.C. and Fracastoro, G. and Proskurnikov, A.V.},
    title={Optimal Clearing Payments in a Financial Contagion Model},
    year={2021},
    journal={Submitted},
    note={online as arXiv:2103.10872}
}

@book{HararyBook:1965,
    author={F. Harary and R.Z. Norman and D. Cartwright},
    title={Structural models. An introduction to the theory of directed Graphs},
    year={1965},
    publisher={Wiley \& Sons},
    address={New York, London, Sydney},
}

@article{EisNoe:01,
  title={Systemic risk in financial systems},
  author={Eisenberg, Larry and Noe, Thomas H},
  journal={Management Science},
  volume={47},
  number={2},
  pages={236--249},
  year={2001},
  publisher={INFORMS}
}

@article{suzuki2002valuing,
  title={Valuing corporate debt: the effect of cross-holdings of stock and debt},
  author={Suzuki, Teruyoshi},
  journal={Journal of the Operations Research Society of Japan},
  volume={45},
  number={2},
  pages={123--144},
  year={2002},
  publisher={The Operations Research Society of Japan}
}

@techreport{alvarez2015mandatory,
  title={Mandatory disclosure and financial contagion},
  author={Alvarez, Fernando and Barlevy, Gadi},
  year={2015},
  institution={National Bureau of Economic Research}
}

@article{sonin2017banks,
  title={Banks as tanks: A continuous-time model of financial clearing},
  author={Sonin, Isaac M and Sonin, Konstantin},
  journal={arXiv preprint arXiv:1705.05943},
  year={2017}
}

@article{feinstein2021dynamic,
  title={Dynamic default contagion in heterogeneous interbank systems},
  author={Feinstein, Zachary and S{\o}jmark, Andreas},
  journal={SIAM Journal on Financial Mathematics},
  volume={12},
  number={4},
  pages={SC83--SC97},
  year={2021},
  publisher={SIAM}
}

@article{capponi2015systemic,
  title={Systemic risk mitigation in financial networks},
  author={Capponi, Agostino and Chen, Peng-Chu},
  journal={Journal of Economic Dynamics and Control},
  volume={58},
  pages={152--166},
  year={2015},
  publisher={Elsevier}
}

@article{ferrara2019systemic,
  title={Systemic illiquidity in the interbank network},
  author={Ferrara, Gerardo and Langfield, Sam and Liu, Zijun and Ota, Tomohiro},
  journal={Quantitative Finance},
  volume={19},
  number={11},
  pages={1779--1795},
  year={2019},
  publisher={Taylor \& Francis}
}

@article{kusnetsov2019interbank,
  title={Interbank clearing in financial networks with multiple maturities},
  author={Kusnetsov, Michael and Maria Veraart, Luitgard Anna},
  journal={SIAM Journal on Financial Mathematics},
  volume={10},
  number={1},
  pages={37--67},
  year={2019},
  publisher={SIAM}
}

@article{chen2021financial,
  title={Financial Network and Systemic Risk—A Dynamic Model},
  author={Chen, Hong and Wang, Tan and Yao, David D},
  journal={Production and Operations Management},
  volume={30},
  number={8},
  pages={2441--2466},
  year={2021},
  publisher={Wiley Online Library}
}

@article{feinstein2020capital,
  title={Capital regulation under price impacts and dynamic financial contagion},
  author={Feinstein, Zachary},
  journal={European Journal of Operational Research},
  volume={281},
  number={2},
  pages={449--463},
  year={2020},
  publisher={Elsevier}
}

@article{banerjee2018dynamic,
  title={Dynamic clearing and contagion in financial networks},
  author={Banerjee, Tathagata and Bernstein, Alex and Feinstein, Zachary},
  journal={arXiv preprint arXiv:1801.02091},
  year={2018},
}

@article{cifuentes2005liquidity,
  title={Liquidity risk and contagion},
  author={Cifuentes, Rodrigo and Ferrucci, Gianluigi and Shin, Hyun Song},
  journal={Journal of the European Economic Association},
  volume={3},
  number={2-3},
  pages={556--566},
  year={2005},
  publisher={Oxford University Press}
}

@article{shin2008risk,
  title={Risk and liquidity in a system context},
  author={Shin, Hyun Song},
  journal={Journal of Financial Intermediation},
  volume={17},
  number={3},
  pages={315--329},
  year={2008},
  publisher={Elsevier}
}

@book{elsinger2009financial,
  title={Financial networks, cross holdings, and limited liability},
  author={Elsinger, Helmut and others},
  year={2009},
  publisher={Oesterreichische Nationalbank Austria}
}

@article{rogers2013failure,
  title={Failure and rescue in an interbank network},
  author={Rogers, Leonard CG and Veraart, Luitgard AM},
  journal={Management Science},
  volume={59},
  number={4},
  pages={882--898},
  year={2013},
  publisher={INFORMS}
}

@article{fischer2014no,
  title={No-arbitrage pricing under systemic risk: Accounting for cross-ownership},
  author={Fischer, Tom},
  journal={Mathematical Finance: An International Journal of Mathematics, Statistics and Financial Economics},
  volume={24},
  number={1},
  pages={97--124},
  year={2014},
  publisher={Wiley Online Library}
}

@article{GlasYou:16,
  title={Contagion in financial networks},
  author={Glasserman, Paul and Young, H Peyton},
  journal={Journal of Economic Literature},
  volume={54},
  number={3},
  pages={779--831},
  year={2016}
}

@article{amini2016fully,
  title={To fully net or not to net: Adverse effects of partial multilateral netting},
  author={Amini, Hamed and Filipovi{\'c}, Damir and Minca, Andreea},
  journal={Operations Research},
  volume={64},
  number={5},
  pages={1135--1142},
  year={2016},
  publisher={Informs}
}

@article{csoka2018decentralized,
  title={Decentralized clearing in financial networks},
  author={Cs{\'o}ka, P{\'e}ter and Jean-Jacques Herings, P},
  journal={Management Science},
  volume={64},
  number={10},
  pages={4681--4699},
  year={2018},
  publisher={INFORMS}
}

@article{Massai:2021,
    author = {Leonardo Massai and Giacomo Como and Fabio Fagnani},
    title = {Equilibria and Systemic Risk in Saturated Networks},
    journal = {Mathematics of Operation Research},
    year = {2021},
    doi={10.1287/moor.2021.1188},
    note={published online, as arXiv:1912.04815}
}

@book{BermanPlemmons_Book,
    author={Berman, Abraham and Plemmons, Robert J.},
    title={Nonnegative Matrices in the Mathematical Sciences},
    publisher={SIAM},
    year={1994},
    address={Philadelphia,PA},
}
\end{small}

\appendix
\section{Appendix: proofs}

\subsection{Technical preliminaries}

We start with some auxiliary material, which will be used in the proofs. 

The following proposition follows, e.g., from~\cite[Corollary~4.3a']{HararyBook:1965}.
\begin{proposition}\label{prop.sink}
Each graph contains at least one sink component. Any strong component that is not a sink is connected to at least one of the sink components by a path.
\end{proposition}

We will also employ several technical propositions, dealing with substochastic\footnote{A \emph{nonnegative} square matrix $A\in\Real{\calV\times\calV}$ is \emph{substochastic} if $\sum_{j\in\calV} a_{ij}\leq 1\,\forall i\in\calV$.} matrices.  
\begin{proposition}\label{prop.submatr-sink}
Let $A=(a_{ij})_{i,j\in\calV}$ be a stochastic matrix and $\calV'\subsetneq \calV$. Then, submatrix
$A'=(a_{ij})_{i,j\in\calV'}$ is \textbf{not} Schur stable if and only if $\calV'$ contains all nodes of some strongly connected sink component of $\calG[A]$.
\end{proposition}
\begin{pf}
The ``if'' part is obvious. If $\calV^0\subseteq\calV'$ is the set of nodes of some sink component, then
$A^0=(a_{i,j})_{i,j\in \calV^0}$ is a stochastic matrix and $A$ is decomposed as
\beq\label{eq.decompose}
A=\begin{pmatrix}
A^0 & \mathbf{O}\\
* & *
\end{pmatrix},
\eeq
where $\mathbf{O}$ is the block of zeros and symbols $*$ denote some submatrices.
Therefore, $A$ has eigenvalue $1$ and is not Schur stable.
The ``only if'' part is implied by~\citep[Lemma~6]{Journal2021}. Thanks to this lemma, $A'$ is not Schur stable if and only if $\calV^0\subseteq\calV'$ exists such that $A^0=(a_{i,j})_{i,j\in \calV^0}$
is a stochastic matrix, which also implies that $A$ is decomposed as in~\eqref{eq.decompose}. In other words, the set of nodes $\calV^0\subseteq\calV'$ is ``closed'': each arc of $\calG[A]$ starting in $\calV^0$ ends also in $\calV^0$. Hence, strong components of graph $\calG[A^0]$ are also strong components of $\calG[A]$, and (due to Proposition~\ref{prop.sink}) at least one of them is a sink.\qed
\end{pf}
\begin{proposition}\label{prop.eq}
Suppose that a substochastic matrix
$A=(a_{ij})_{i,j\in\calV'\cup\calV''}$ is Schur stable, where
$\calV'\cap\calV''=\emptyset$. Then a vector $\xi\geq 0$ exists
such that
\[
\xi_i-(A\tran\xi)_i
\begin{cases}
>0,\,i\in\calV',\\
=0,\,i\in\calV''.
\end{cases}
\]
\end{proposition}
\begin{pf}
If $A$ is Schur stable, then $(I-A)^{-1}=\sum_{k=0}^{\infty}A^k\geq 0$ exists.
Choosing an arbitrary vector $e$ such that $e_i>0\,\forall
i\in\calV'$ and $e_i=0\,\forall i\in\calV''$, the vector
$\xi=(I-A\tran)^{-1}e$ thus is also nonnegative. By construction,
$(\xi-A\tran\xi)_i=e_i>0$ if and only if $i\in\calV'$.\qed
\end{pf}
We also need a special form of the Perron-Frobenius theorem.
\begin{lemma}\label{lem.pf}
Let $A\in\Real{\calV\times\calV}$ be a stochastic matrix and $\calV^0$ be the set of nodes of some sink component in $\calG[A]$. Then, vector $\pi\in\Real{\calV}$ exists such that
\beq\label{eq.eigenvec}
A^{\top}\pi=\pi,\;\pi\tran\one=1,\;\pi_i
\begin{cases}
>0,\,\forall i\in\calV^0,\\
=0,\,\forall i\not\in\calV^0.
\end{cases}
\eeq
\end{lemma}
\begin{pf}
The special case where $\calG[A]$ is a strongly connected graph ($A$ is irreducible) and $\calV^0=\calV$ is immediate from Perron-Frobenius theorem for irreducible matrices~\cite[Ch.2, Th.~1.3]{BermanPlemmons_Book}.
Otherwise, matrix $A$ has structure~\eqref{eq.decompose}, where $A^0=(a_{ij})_{i,j\in\calV^0}$
is irreducible (graph $\calG[A^0]$ is strongly connected by the definition of a strongly connected component). Introducing the Perron-Frobenius eigenvector $\pi^0>0$ of matrix $A^0$, vector $\pi$ can be defined as follows:
\[
\pi_i=\pi_i^0\,\forall i\in\calV^0,\quad \pi_i=0\,\forall i\not\in\calV^0.
\]
\end{pf}

Lemmas~\ref{lem.unique} and~\ref{lem.pf} have a simple corollary, which will be used in the proof of Theorem~\ref{thm.prorata}.
\begin{corollary}\label{cor.max}
Given a stochastic matrix $A$ and nonnegative vectors $c,\bar p\geq 0$, consider the maximal payment vector $p^*=p^*[A,c,\bar p]$ from Lemma~\ref{lem.unique}. Suppose that graph $\calG[A]$ has a strongly connected component with set of nodes $\calV^0\subset\calV$, which is a sink (no arc leaves it) and is such that $\bar p_i>0\,i\in\calV^0$. for Then $p_i^*>0\,\forall i\in\calV^0$.
\end{corollary}
\begin{pf}
Choosing $\pi$ as in~\eqref{eq.eigenvec}, one has $\ve\pi\in\calP\ped{pr}(c,\bar p)$ for $\ve>0$ small enough
(so small that $\ve\pi_i<\bar p_i\,\forall i\in\calV^0$). Since $p^*$ is the maximal element of $\calP\ped{pr}(c,\bar p)$, we have $p_i^*\geq\ve\pi_i>0\,\forall i\in\calV^0$.\qed
\end{pf}

\subsection{Proof of Theorem~\ref{thm.1}} \label{sec.proofs_th1}
We introduce the following notation: for a pair of banks $i$, $j\ne i$ let
\[
\delta_{ij}(t)=\bar p_{ij}(t)-p_{ij}(t)\geq 0,\;t\in \calT
\]
be the amount bank $i$ owes to bank $j$ before period $t+1$. In view of~\eqref{eq_pdynamics}, $\bar p_{ij}(t+1)>0$ if and only if $\delta_{ij}(t)>0$.

The proof is based on a simple transformation, which we call the transformation of
advance payment (TAP). Let $J$ be a subset of arcs in graph $\calG[P(t_0)]$, where $1\leq t_0\leq T$, and $\ve>0$.
For $(i,j)\in J$, one has $\bar p_{ij}(t_0)\geq p_{ij}(t_0)>0$ and, thus $\delta_{ij}(t_0-1)>0$.
The TAP with parameters $(t_0,\ve,J)$ modifies matrices $P(t_0-1)$ and $P(t_0)$ as follows:
\begin{itemize}
\item at time $t_0-1$, payment on each arc from $J$ is increased
\[
p_{ij}(t_0-1)\mapsto 
p_{ij}(t_0-1)+
\alpha^{-1}\varepsilon\;\;\forall (i,j)\in J;
\]
\item at time $t_0$, payment on each arc from $J$ is decreased
\[
p_{ij}(t_0)\mapsto 
p_{ij}(t_0)-\varepsilon\;\;\forall (i,j)\in J;
\]
\item  all other entries of $P(t_0-1)$ and $P(t_0)$ and remaining matrices $P(t)$, $t\ne t_0-1,t_0$ remain unchanged.
\end{itemize}
Obviously, this transformation increases the objective function~\eqref{eq:prob01a} by $(\alpha^{-1}a_{t_0-1}- a_{t_0})|J|>0$. For $\ve>0$ being sufficiently small, the TAP transformation preserves constraints~\eqref{eq.cond-p-1+}: it suffices to choose $\ve<\min\{p_{ij}(t_0):(i,j)\in J\}$. Conditions~\eqref{eq.cond-p-1a} also retain their validity, provided that $\ve<\alpha\min\{\delta_{ij}(t_0-1):(i,j)\in J\})$.
Notice that the nominal payment matrices $\bar P(0),\ldots\bar P(t_0)$ remain unchanged, and hence the condition $P(t)\leq\bar P(t)$ (equivalent to~\eqref{eq.cond-p-1a}) holds for all $t\leq t_0$. For $t\geq t_0$, the sum in the left-hand side of~\eqref{eq.cond-p-1a} is invariant under the TAP transformation, so the constraint is also not violated.
Finally, constraints~\eqref{eq.cond-p-2a} (equivalent to $w_i(t+1)\geq 0$) also hold for all $t$ except for, possibly, $t=t_0-1$ and $t=t_0$, because other matrices $P(t)$ remain unchanged.

In view of the optimality of sequence $[P]$, the TAP transformation with parameters $(t_0,J,\ve)$, where $\ve>0$ is sufficiently small, violate~\eqref{eq.cond-p-2a} at $t=t_0-1$ or at $t=t_0$.

\textbf{Step 1.} We first prove the last statement of Theorem~\ref{thm.1}. Assume that this statement is not valid and
a cycle $i_1\xrightarrow[]{}\ldots \xrightarrow[]{} i_s \xrightarrow[]{} i_1$ exists in $\calG[P(t_0)]$, where $t_0\geq 1$. Choosing the set of arcs $J=\{i_1\xrightarrow[]{} i_2,\ldots,i_{s-1}\xrightarrow[]{}i_s, i_s\xrightarrow[]{}i_1\}$,
the TAP transformation with parameters $(t_0,J,\ve)$ (with $\ve>0$ small enough), obviously, leaves
the vectors $\phi\ap{out}(t)-\phi\ap{in}(t)=P(t)\one-P(t)^{\top}\one$ unchanged, and thus constraints~\eqref{eq.cond-p-2a} are not violated, which
leads one to a contradiction with the optimality of $[P]$.

\textbf{Step 2.} Suppose now that $t_*=t_*(i)\leq T$ is defined as described in Theorem~\ref{thm.1}  yet~\eqref{eq.default} fails to hold at some period $0\leq t<t_*$.
Let $t_+<t_*$ be the \emph{last} period when~\eqref{eq.default} fails, that is, the maximum of $t<t_*$ such that
$w_i(t+1)=\phi_i\ap{in}(t)+w_i(t)-\phi_i\ap{out}(t)>0$.

Notice first that $t_+<T-1$. Otherwise, one would have $w_i(T)>0$ and $t_*=T$, in particular, $\delta_{ij}(T-1)>0$ for some $j\ne i$. Increasing $p_{ij}(T-1)$ by a sufficiently small value $\ve>0$, one could obviously preserve all constraints and also increase the objective function.

Denoting for brevity $t_0\doteq 1+t_+\leq t_*$, one thus has $t_0<T$. The definition of $t_+$ and $t_0$ implies that $\phi_i\ap{out}(t_0)>0$. Indeed, if $t_0=t_*$,
then one has $\phi_i\ap{out}(t_0)=\bar \phi_i\ap{out}(t_*)>0$ by definition of $t_*$. Otherwise,
$\phi_i\ap{out}(t_0)=\phi_i\ap{in}(t_0)+w_i(t_0)\geq w_i(t_0)=w_i(1+t_+)>0$ due to the choice of $t_0$.

We know that graph $\calG[P(t_0)]$ contains no cycles and, in particular, all its strongly connected components are trivial (single-node) graphs. Since $\phi_i\ap{out}(t_0)>0$, node $i$ is not a sink node.
Proposition~\ref{prop.sink} ensures that $i$ is connected to a sink node $k$ by a path
$i\xrightarrow[]{}j_1\xrightarrow[]{}\ldots\xrightarrow[]{} j_s\xrightarrow[]{}k$
(all nodes $i,j_1,\ldots,j_s,k$, $s\geq 0$ are mutually different).
Let $J$ be the set of arcs in this path. The TAP with parameters
$(t_0,J,\ve)$ with $\ve>0$ small enough, obviously, preserves~\eqref{eq.cond-p-2a} (equivalent to $w(t+1)\geq 0$) at $t=t_0-1$ or at $t=t_0$.
Indeed, the TAP leaves the components $w_j(t_0),w_j(t_0+1)$ for each $j\ne i,k$ invariant.
The component $w_k(t_0)$ increases (becoming thus positive), and hence $w_k(t_0+1)$ is also positive (recall that $\phi_k\ap{out}(t_0)=0$). The TAP transformation decreases $w_i(t_0)$ by $\alpha^{-1}\varepsilon$ (providing that $w_i(t_0)>0$ for $\ve$ being small), however, $\phi_i\ap{out}({t_0})$ is decreased by $\ve$, so that $w_i(t_0+1)$ is increased by $(1-\alpha^{-1})\varepsilon\geq 0$, and inequality $w_i(t_0+1)\geq 0$ is preserved. Hence,  constraints~\eqref{eq.cond-p-2a} are not violated, and we arrive at a contradiction with optimality of $[P]$.

\textbf{Step 3.} The proof of implication~\eqref{eq.default-imp} is now straightforward. Suppose that $\phi_i\ap{out}(t)<\bar\phi_i\ap{out}(t)$ at some period $t\in\mathcal{T}$. Then, obviously, $t<t_*(i)$, and hence
$\phi_i\ap{in}(t)+w_i(t)-\phi_i\ap{out}(t)=0$ due to~\eqref{eq.default}.
\qed

\subsection{Proof of Theorem~\ref{thm.prorata}}\label{sec.proofs}


In the proof, we will use a \emph{transformation of advanced payment} (TAP), which is similar to the transformation used in the proof of Theorem~\ref{thm.1}. The TAP is determined by time instant $t_*$, scalar $\ve>0$ and non-negative vector $\zeta\geq 0$; it
replaces sequence $[p]$ by the sequence
$[\hat p]$, where
\beq\label{eq.tap}
\hat p(t)=
\begin{cases}
p(t),\,&t\ne t_*,t_*+1,\\
p(t_*)+\ve\alpha^{-1}\zeta,\,&t=t_*,\\
p(t_{+})-\ve\zeta,\,&t=t_*+1.
\end{cases}
\eeq
In other words, some payments are transferred (taking into account the interest rate $\alpha\geq 1$) from period $t_*+1$
to the previous period $t_*$.

If $[p]$ satisfies constraints~\eqref{eq.cond-p-1++}-\eqref{eq.cond-p-2b}, then $[\hat p]$ also obeys all constraints, except for, possibly:
1) constraint~\eqref{eq.cond-p-1++} at $t_*+1$ (at other periods, $\hat p(t)\geq p(t)$);
2) constraint~\eqref{eq.cond-p-1b} at $t=t_*$ (at other periods, the left-hand side of~\eqref{eq.cond-p-1b} remains invariant under the TAP);
3) constraints~\eqref{eq.cond-p-2b} at periods $t=t_*,\ldots,T-1$ (for $t<t_*$, the left-hand side of~\eqref{eq.cond-p-2b} remains invariant under the TAP).
 Also, for any $\zeta\neq 0$ and $\ve>0$ the TAP always increases the value of the objective function~\eqref{eq:prob01b}, because $a_{t_*}>\alpha a_{t_*+1}$.

For the optimal sequence of payment vectors $[p^*]$, we are going to prove that $p^*(t)$ (at each $t$) is a maximizer at problem~\eqref{eq:clearingopt-b},\eqref{eq:clearingopt-c}, or, equivalently,~\eqref{eq.p*-dynamic} holds,
via backward induction on $t=T-1,T-2,\ldots,0$. Here $w^*(t)$ is the net worth~\eqref{eq.w-evolves} corresponding to $p^*(t)$.

\textbf{The induction base} $t=T-1$ is obvious, recalling that constraints~\eqref{eq:clearingopt-c} are equivalent to~\eqref{eq.cond-p-1b},\eqref{eq.cond-p-2b}. If $p^*(T-1)$ were not a maximizer in~\eqref{eq:clearingopt-b},\eqref{eq:clearingopt-c} with $T-1$, the value of objective function in~\eqref{eq:prob01b} could be increased.

\textbf{The induction step.} Suppose that our statement has been proved for $t=t_*+1,\ldots,T-1$. In particular, at each $t>t_*$ vector $p^*(t)$ obeys the equation~\eqref{eq:dynfixpoint}.
We are now going to prove that~\eqref{eq.p*-dynamic} holds at $t=t_*$.
The proof is based on Lemma~\ref{lem.unique} and is performed in two steps.

\textbf{Step 1.} We first show that  each strongly connected sink component of $\calG[A]$ contains node $i$ such that $p_i^*(t_*)=\bar p_i^*(t_*)$. Suppose that the statement is not correct and consider
such a sink strong component of $\calG[A]$ with the set of nodes $\calV^0\subseteq\calV$ that
$\bar p_i^*(t_*)>p_i^*(t_*)\,\forall i\in\calV^0$, or, equivalently, $\bar p_i^*(t_*+1)>0\,\forall i\in\calV^0$. Applying Corollary~\ref{cor.max} to $\bar p=\bar p(t_*+1)$ and recalling that~\eqref{eq.p*-dynamic} holds at $t=t*+1$, one has $p_i(t+1)>0\,\forall i\in\calV^0$.

Introducing the eigenvector from Lemma~\ref{lem.pf}, consider the TAP~\eqref{eq.tap} with $p=p^*$, $\zeta=\pi$ and $\ve>0$ sufficiently small. Since $\zeta=A\tran\zeta$,
the left-hand side of~\eqref{eq.cond-p-2b} remains invariant under the TAP, and hence $[\hat p]$ obeys constraints~\eqref{eq.cond-p-2b}.
Since $\zeta_i=\pi_i=0$ for $i\not\in\calV^0$ and $p_i(t+1)>0\,\forall i\in\calV^0$, constraint~\eqref{eq.cond-p-1++} at $t=t_*+1$ is also preserved by the TAP when $\ve>0$ is so small that
$\hat p_i(t_*+1)=p_i^*(t_*+1)-\ve\zeta_i>0\,\forall i\in\calV^0$. Finally, constraint~\eqref{eq.cond-p-1b} at $t=t_*$ can be rewritten as $p_i(t_*)\leq\bar p_i(t_*)$. Recalling that
$\zeta_i=\pi_i=0$ for $i\not\in\calV^0$ and $p_i^*(t_*)<\bar p_i^*(t_*)$ for $i\in\calV^0$, it is obvious that the TAP does not violate this constraint for $\ve>0$ sufficiently small.
As has been noticed, the remaining constraints are always preserved by the TAP.
The new sequence of payment vectors $[\hat p]$ thus satisfies all the constraints~\eqref{eq.cond-p-1++}-\eqref{eq.cond-p-2b}
and corresponds to a \emph{larger} value of the objective function, which leads to a contradiction with the optimality of $[p^*]$. The contradiction shows that inequality $p_i^*(t_*)<\bar p_i^*(t_*)$ is violated for at least one index $i\in\calV^0$.

\textbf{Step 2.} In view of Lemma~\ref{lem.unique}, statement (d) (applied for $\bar p=\bar p^*(t_*)$ and $c=c(t_*)+w(t_*)$), to prove~\eqref{eq.p*-dynamic} at $t=t_*$ it remains to prove~\eqref{eq:dynfixpoint} at $t=t_*$.

Assume that~\eqref{eq:dynfixpoint} fails to hold, that is, index $s\in\calV$ exists such that $p_s^*(t_*)<\bar p_s^*(t_*)$ and $p_s^*(t_*)<c_s(t_*)+w_s^*(t_*)+(A\tran p^*(t_*))_s$.
We are going to show that this leads to a contradiction with the assumption that sequence $[p^*]$ is optimal, using the TAP~\eqref{eq.tap}.

We first define the following sets of indices. Let $\calV^0\ne\emptyset$ consist of such nodes $i$ that $p_i^*(t_*)=\bar p_i^*(t_*)$ (at Step 1, we have shown every strongly connected sink component of $\calG[A]$ contains an element from $\calV^0$) and $\tilde\calV\doteq\calV\setminus\calV^0$. Obviously, $s\in\tilde\calV$. We introduce the submatrix $\tilde A=(a_{ij})_{i,j\in\tilde\calV}$ and the corresponding graph $\tilde\calG=\calG[\tilde A]$. Let $\calV^1$ stand for all nodes $i\in\tilde\calV$, $i\ne k$ that $i$ \emph{are not} reachable from $s$ in $\tilde G$,
$\calV^2$ stand for all nodes $i\in\tilde\calV$, $i\ne s$ that $i$ \emph{are} reachable from $s$ in $\tilde G$.

\textbf{Step 2a.}  We first show that $p_i^*(t_*+1)>0$ for all $i\in\calV^2$.

By construction, $\bar p_i^*(t_*+1)=\alpha[\bar p_i^*(t_*)-p_i^*(t_*)]>0\,\forall i\in\tilde\calV$. The induction hypothesis entails now that $p_i^*(t_*+1)>0$ for $i\in\calV^2\cup\{s\}$.
Indeed, by assumption $w_s^*(t_*+1)=c_s(t_*)+w_s^*(t_*)+(A\tran p^*(t_*))_s-p_s^*(t_*)>0$. Recalling that~\eqref{eq:dynfixpoint} holds at $t=t_*+1$, one shows that $p_s^*(t_*+1)>0$.
If node $\ell\in\calV^2$ is directly accessible from $s$ (that is, $a_{sl}>0$) in $\tilde\calG$, then~\eqref{eq:dynfixpoint} at $t=t_*+1$ implies that $p_{\ell}^*(t_*+1)>0$, because
$(A^{\top}p^*)_{\ell}\geq a_{s\ell}p_s^*(t_*+1)>0$. Similarly, if a path $s\xrightarrow[]{}\ell\xrightarrow[]{}m$ exists in $\tilde\calG$, then $p_{m}^*(t_*+1)>0$ due to~\eqref{eq:dynfixpoint}, because $a_{\ell m}p_{\ell}^*(t_*+1)>0$, and so on: via induction of the length of the path connecting $s$ to $i\in\calV^2$, one shows that $p_i^*(t_*+1)>0$ for all $i\in\calV^2$.

\textbf{Step 2b.} As has been shown at Step~1, set $\tilde\calV$  does not contain any strongly connected sink component of $\calG[A]$, and hence matrix $\tilde A$ and all its submatrices are Schur stable (Proposition~\ref{prop.submatr-sink}). Applying Proposition~\eqref{prop.eq} to $\calV'=\{s\}$ and $\calV''=\calV^2$, a vector $\xi\in\Real{\calV^2\cup\{s\}}$ exists such that
\beq\label{eq.aux2}
\xi\geq 0\;\;\text{and}\;\;\xi_i-\sum_{j\in\calV^2\cup\{s\}}a_{ji}\xi_j
\begin{cases}
>0,\,i=s;\\
=0,\,i\in\calV^2.
\end{cases}
\eeq
Define the vector $\zeta$ as follows: $\zeta_i\doteq 0$ for $i\in\calV^0\cup\calV^1$ and $\zeta_i\doteq\xi_i$ for $i\in\calV^2\cup\{s\}$.
Then,
\beq\label{eq.aux2+}
\begin{gathered}
\zeta_s-(A\tran\zeta)_s>0\\
\zeta_i-(A\tran\zeta)_i=0\quad\forall i\in\calV^1\cup\calV^2.
\end{gathered}
\eeq
Indeed, for $i\in\calV^2\cup\{s\}$ one has
\[
(A\tran\zeta)_i=\sum_{j\in\calV^2\cup\{s\}}a_{ji}\zeta_j+\sum_{j\in\calV^0\cup\calV^1}a_{ji}\underbrace{\zeta_j}_{=0}\leq \zeta_i
\]
due to~\eqref{eq.aux2}, which inequality can be strict only when $i=s$.  Obviously, if $i\in\calV^1$ and $j\in\calV^2\cup\{s\}$, then $a_{ji}=0$ (otherwise, $i$ would be reachable from $s$ in graph $\tilde G$, contradiction to the definition of $\calV^1$). Thus,
\[
(A\tran\zeta)_i=\sum_{j\in\calV^2\cup\{s\}}\underbrace{a_{ji}}_{=0}\zeta_j+\sum_{j\in\calV^0\cup\calV^1}a_{ji}\underbrace{\zeta_j}_{=0}=0=\zeta_i.
\]
for all $i\in\calV^1$.

\textbf{Step 2c.} We are now ready to show that the transformation~\eqref{eq.tap} with the constructed vector $\xi$ and $p=p^*$ does not violate constraints~\eqref{eq.cond-p-1++}-\eqref{eq.cond-p-2b} if $\ve>0$ is small.

As we know, constraint~\eqref{eq.cond-p-1++} has to be checked only at $t=t_*+1$. By construction, $\zeta_i=0$ unless $i\in\calV^2\cup\{s\}$.
As we have seen at Step 2a, for such indices one has $p_i^*(t_*+1)>0$. Hence,
$\hat p_i(t_*+1)=p_i^*(t_*+1)-\ve\zeta_i\geq 0\,\forall i\in\calV$ for $\ve>0$ being sufficiently small.

Constraint~\eqref{eq.cond-p-1b} has to be tested only at $t=t_*$. For each $i\in\calV^2\cup\{s\}$ one has $0<\bar p_i^*(t_*)-p_i^*(t_*)=\alpha^{t_*} \bar p_i-\sum\nolimits_{k=0}^{t_*} \alpha^{t_*-k}p_i^*(k)$, which inequality, obviously, remains valid also when $p^*$ is replaced by $[\hat p]$ (provided that $\ve>0$ is small enough).
On the other hand, for $i\in\calV^0\cup\calV$, \[
0\leq\alpha^{t_*} \bar p_i-\sum_{k=0}^{t_*} \alpha^{t_*-k}p_i^*(k)=\alpha^{t_*} \bar p_i-\sum_{k=0}^{t_*} \alpha^{t_*-k}\hat p_i^*(k).
\]
Hence, $[\hat p]$ satisfies constraints~\eqref{eq.cond-p-1b} if $\ve>0$ is small.

Finally, we have to check constraints~\eqref{eq.cond-p-2b} for all $t\geq t_*$. Recall that
$w(t+1)=C(t) + \sum\nolimits_{k=0}^{t} \left(A\tran p(k)- p(k) \right)$
due to~\eqref{eq.w-evolves} and~\eqref{eq:proratac}, and~\eqref{eq.cond-p-2b} is equivalent to the inequality $w(t+1)\geq 0$. Denoting the net worth vectors corresponding to $[\hat p]$ by
\beq\label{eq.aux3}
\begin{split}
\hat w(t+1)\doteq C(t) + \sum\nolimits_{k=0}^{t} \left(A\tran \hat p(k)- \hat p(k) \right)=\\
=\hat w(t)+c(t)+A\tran \hat p(t)- \hat p(t),
\end{split}
\eeq
our goal is to show that $\hat w(t+1)\geq 0$ for $t=t_*,\ldots,T-1$ and $\ve>0$ being small.

By assumption, for each $i\in\calV^0$ and each $t>t_*$ one has $\hat p_i^*(t)=p_i^*(t)=0$ (node $i$ pays full debt at period $t_*$). As has been already shown, we have $\hat p(t)\geq 0$ for all $t\geq 0$.  In view of~\eqref{eq.aux3}, for every such index the sequence $\hat w_i(t)$ is non-decreasing as $t=t_*,t_*+1,\ldots,T-1$:
\[
\begin{gathered}
\hat w_i(t+1)=\hat w_i(t)+c_i(t)+(A\tran \hat p(t))_i-\hat p_i(t)\geq \hat w_i(t)\\
\forall t\geq t_*+1,\,\forall i\in\calV^0.
\end{gathered}
\]
On the other hand, $\hat p(t_*)\geq p^*(t_*)$ and $\hat p_i(t_*)=p_i^*(t_*)\,\forall i\in\calV^0$, whereas $p^*(t)=\hat p(t)$ for $t<t_*$. In view of this, $\hat w_i(t_*+1)=w_i^*(t_*+1)\geq 0\,\forall i\in\calV^0$, which shows that $\hat w_i(t+1)\geq 0$ for $t=t_*,\ldots, (T-1)$ and $i\in\calV^0$.

On the other hand,~\eqref{eq.tap} entails that
\[
\hat w(t+1)-w^*(t+1)=
\begin{cases}
0,\,t<t_*,\\
-\alpha^{-1}\ve(\zeta-A\zeta),\,t=t_*\\
(1-\alpha^{-1})\ve(\zeta-A\zeta),t>t_*.
\end{cases}
\]
Since $1-\alpha^{-1}\geq 0$, inequalities~\eqref{eq.aux2+} entails that $\hat w_i(t+1)\geq w_i^*(t+1)\geq 0$ for $i\in\calV^1\cup\calV^2\cup\{k\}$ and $t>t_*$.
Furthermore, $\hat w_i(t_*+1)=w_i^*(t_*+1)\geq 0$ for $i\in\calV^1\cup\calV^2$. Finally, by assumption $w_s^*(t_*+1)>0$ entails that $\hat w_s(t_*+1)>0$ provided that $\ve>0$ is small enough.
We have demonstrated that $\hat w(t+1)\geq 0$ (equivalently, $[\hat p]$ satisfies~\eqref{eq.cond-p-2b}) for $t=t_*,\ldots, (T-1)$, provided that $\ve>0$ is chosen sufficiently small.

The assumption about the existence of index $s$ such that $p_s^*(t_*)<\bar p_s^*(t_*)$ and $p_s^*(t_*)<c_s(t_*)+w_s^*(t_*)+(A\tran p^*(t_*))_s$ has led us to the contradiction with optimality of sequence $[p]$,
because the new sequence $[\hat p]$ satisfies all constraints and corresponds to a large value of the objective function. Therefore,~\eqref{eq:dynfixpoint} should hold at $t=t_*$, which,
along with the statement proved at Step~1 and Lemma~\ref{lem.unique}, ensures that~\eqref{eq.p*-dynamic} holds at $t=t_*$. This finishes the proof of induction step.

The uniqueness of the optimal solution is now trivial. Lemma~\ref{lem.unique} ensures the uniqueness of $p^*(0)$, which is the  maximizer at~\eqref{eq:clearingopt-b},\eqref{eq:clearingopt-c} with $t=0$
(and depends only on $c(0)$). Similarly, $p^*(1)$ (depending on $c(1)$ and $p^*(0)$) is uniquely found as the maximizer at~\eqref{eq:clearingopt-b},\eqref{eq:clearingopt-c} with $t=1$, and so on; using induction on
$t=0,\ldots,T-1$, one shows that $p^*(t)$ is defined uniquely and depends on $c(t)$ and $p^*(0),\ldots,p^*(t-1)$.
\qed

\subsection{A remark on the structure of the cost function}

Note that the proofs in the previous subsections do not use the representation of coefficients $a_t$ in~\eqref{eq:prob01a} and~\eqref{eq:prob01b}.
The constants~\eqref{eq:ai} can be replaced by any positive numbers $a_0,\ldots,a_{T-1}$ such that $a_{t-1}>\alpha a_{t}\,\forall t=1,\ldots,T-1$.

In particular, instead of minimizing the loss function~\eqref{eq.cost}, one may minimize a more general function
\begin{gather}
J([P])=(1-\eta)L([P])+\eta\one^{\top}\bar P(T)\one,\label{eq.cost-eta}\\
\bar P(T)\overset{\eqref{eq.bar-p-evolves}}{=}\alpha^T \bar P -  \sum\nolimits_{k=0}^{T-1} \alpha^{T-k}P(k),\notag
\end{gather}
where $\eta\in[0,1)$. The loss function~\eqref{eq.cost} corresponds to $\eta=0$; the weight $\eta>0$
corresponds to the additional penalty on unpaid liabilities (recall that $\bar P(T)=0$ if and only if there is no default at the terminal time).
To minimize the cost function~\eqref{eq.cost-eta}, one has to maximize the function~\eqref{eq:prob01a} (or, in the case of pro-rata payments,~\eqref{eq:prob01b}) with the weights
\beq\label{eq.ai-eta}
a_t=\eta\alpha^{T-t}+(1-\eta)\sum_{k=0}^{T-t-1}\alpha^k.
\eeq
Theorems~\eqref{thm.1} and~\eqref{thm.prorata} retain their validity, replacing the coefficients~\eqref{eq:ai} by~\eqref{eq.ai-eta}.
\end{document}